\documentclass[10pt]{amsart}
\usepackage{rotating}
\usepackage{mathcomp,amscd}
\usepackage{amssymb}
\usepackage{mathtools}
\usepackage[all]{xy}
\usepackage{subfigure}
\usepackage[margin=1cm]{caption}
\usepackage{wasysym}
 
\usepackage{color}

\newcommand{\Br}{\mathrm{Br}}

\renewcommand{\curve}{\Sigma}
\newcommand{\univ}{\mathrm{univ}}

\newcommand{\pseudosimple}{pseudosimple}
\newcommand{\map}{f}

\newcommand{\ab}{\mathrm{ab}}

\newcommand{\Zprimes}{{\textstyle{\Bbb{Z}}[\frac{1}{\mathcal{P}}]}}

\newcommand{\M}{G}
\newcommand{\MZ}{G'}
\newcommand{\primes}{\mathcal{P}}

\newcommand{\inv}{{\mathrm{inv}}}
\newcommand{\iso}{{\iota}}
\newcommand{\PP}{{\Bbb P}}

\newcommand{\Out}{\mathrm{Out}}

\newcommand{\nuu}{\nu}

\newcommand{\tupleg}{\underline{g}}

\newcommand{\Hur}{{\mbox{\sc Hur}}}
\newcommand{\AHur}{{\sf Hur}}

\newcommand{\Conf}{{\mbox{\sc Conf}}}
\newcommand{\AConf}{{\sf Conf}}

\usepackage[OT2,T1]{fontenc}
\DeclareSymbolFont{cyrletters}{OT2}{wncyr}{m}{n}
\DeclareFontFamily{OT1}{rsfs}{}
\DeclareFontEncoding{OT2}{}{} % to enable usage of cyrillic fonts
  
     \DeclareFontShape{OT1}{rsfs}{n}{it}{<-> rsfs10}{}
\DeclareMathAlphabet{\mathscr}{OT1}{rsfs}{n}{it}

\newcommand{\C}{{\Bbb C}}
\newcommand{\Z}{{\Bbb Z}}

\newcommand{\PGL}{\mathrm{PGL}}

\newcommand{\cG}{\mathcal{G}}

\newcommand{\cF}{\mathcal{F}}
\newcommand{\cH}{\mathcal{H}}
\newcommand{\cI}{\mathcal{I}}
\newcommand{\F}{{\Bbb F}}

\newcommand{\bbP}{{\Bbb P}}
\newcommand{\AP}{{\sf P}}  % sans serif for analytic space; small caps for schemes

\newcommand{\Q}{{\Bbb Q}}
 
\newcommand{\Gal}{\mbox{Gal}}
\newcommand{\Aut}{\mbox{Aut}}
\newcommand{\Inn}{\mbox{Inn}}

\newcommand{\Mab}{G^{\mathrm{ab}}}

\newcommand{\Sym}{\mathrm{Sym}}
\newcommand{\Alt}{\mathrm{Alt}}
\newcommand{\cP}{\mathcal{P}}

\newcommand{\disc}{\mbox{disc}}

\numberwithin{equation}{section}
\numberwithin{table}{section}
\numberwithin{figure}{section}
\newtheorem{Theorem}{Theorem}[section]
\newtheorem*{Theorem*}{Theorem}
\newtheorem{Proposition}[Theorem]{Proposition}
\newtheorem{Corollary}[Theorem]{Corollary}
\newtheorem*{Proposition*}{Proposition}
\newtheorem{Lemma}[Theorem]{Lemma} 
\newtheorem*{lemma*}{Lemma}

\newtheorem{Conjecture}[Theorem]{Conjecture}

\newcommand{\familial}{}
\newcommand{\cmmt}[1]{}
 \allowdisplaybreaks

\swapnumbers
\title{Hurwitz Monodromy and Full Number Fields}

\author{David P.\ Roberts}
\address{Division of Science and Mathematics, University of
  Minnesota Morris, Morris, MN 56267, USA}
\email{roberts@morris.umn.edu}
  
\author{Akshay Venkatesh}
\address{Department of Mathematics, Building 380, Stanford University, Stanford, CA 94305, USA}
\email{akshay@math.stanford.edu}

\begin{document}

\begin{abstract}  
We give conditions for the monodromy group of a
Hurwitz space over the configuration space of branch points to be
 the full alternating or symmetric
group on the degree.   Specializing the resulting coverings suggests the existence
of many number fields with surprisingly little ramification ---
 for example,  the existence of  infinitely many $A_m$ or $S_m$ number fields
unramified away from $\{2,3,5\}$.  
  \end{abstract}
\maketitle

 \tableofcontents

\section{Introduction}
\label{introduction}

\subsection{Overview}  
\label{overview}
       The motivation of this paper is an open problem posed in \cite{MR} 
       concerning number fields, as follows.   
       Say that a degree $m$ number field $K$ is {\em full} if its associated Galois group is either  $A_m$ or $S_m$.  Fix
       a finite set of primes $\cP$. The problem is, {\em Are there infinitely many full fields $K$
       for which the discriminant of $K$ is divisible only by primes in $\cP$?}
       
       In this paper we present a construction, with origins in work of Hurwitz,
       which gives many fields of this type. On the basis of this construction we propose: 
       
       \begin{Conjecture}
     {

\label{mc}
Suppose $\cP$ contains the set of prime divisors of the order of a nonabelian finite simple group.
Then there exist infinitely many full fields unramified outside $\cP$. }
\end{Conjecture}

Our construction amounts to specializing suitable coverings of $\Q$-algebraic varieties
at suitable rational points. 
In the current paper, we analyze the geometric part of the construction, defining
the varieties and proving that the   {\em geometric} monodromy group is $A_m$ or $S_m$.  A sequel paper \cite{HNF} 
provides experimental evidence that fullness is sufficiently preserved by 
the specialization step for Conjecture~\ref{mc} to be true.

       Now some words as to why we find this surprising: 
        In \cite{RWild} the first-named author applied Bhargava's mass heuristic \cite{B}
        to the open question.  
       For given $\cP$,  a finite number was obtained for the total expected number
       of full fields $K$.  Accordingly, it was conjectured in \cite{RWild} that the answer to the 
       question is  {\em no} for all ${\cP}$.   However, the construction given in this paper
      systematically gives fields 
       which escape the influence of the mass heuristic.  
          It is clear from these fields 
       that \cite{RWild} applied the mass heuristic out of its regime 
       of applicability. 
   
   Sections~\ref{general}, \ref{braid} and \ref{lifting} provide short summaries 
 of large theories  and serve to establish our setting.  
 Section~\ref{newtry2} states our main theorem, which we call the
full-monodromy theorem.    It has the form that two statements I and II
are equivalent.  Sections~\ref{direct} and \ref{reverse} prove the theorem by
establishing I $\Rightarrow$ II and II $\Rightarrow$ I respectively.  
 \S\ref{introFM} provides an overview of this material.
 
  \S\ref{introNF} provides an overview of  our  construction of full number fields.   Section~\ref{FullProposal} concludes the paper with more details,
   a sampling of the numerical evidence for Conjecture~\ref{mc}, and 
  further discussion of full number fields in large degree
  ramified within a prescribed $\cP$.

\subsection{The full-monodromy theorem}  
 
 \label{introFM}
 
    Define a {\em Hurwitz
parameter} to be  a triple $h = (G,C,\nu)$ where $G$ is a finite group, $C = (C_1,\dots,C_r)$
is a list of conjugacy classes generating $G$, and $\nu = (\nu_1,\dots,\nu_r)$ is a list of positive integers,
with $\nu$ {\em allowed} in the sense  that $\prod [C_i]^{\nu_i} = 1$ in the abelianization $G^{\rm ab}$.    
A Hurwitz parameter determines
an unramified covering of  complex algebraic
varieties:
  \begin{equation} \label{maincover1} 
  \pi_{h}: \AHur_{h} \rightarrow \AConf_{\nuu}.
  \end{equation}
Here the cover $\AHur_h$ is a Hurwitz variety
parameterizing  certain covers of the complex projective line $\AP^1$, where the coverings are ``of type $h$.'' 
The base $\AConf_{\nuu}$ is the variety whose points are tuples $(D_1,\dots,D_r)$ of disjoint divisors $D_i$ 
 of $\AP^1$, with $\deg(D_i) = \nu_i$.   The map $\pi_h$ sends a cover to its branch 
 locus.   
 
   In complete analogy with the use of the term for number fields, we say that a cover
 of connected complex algebraic varieties $\mathsf{X} \rightarrow \mathsf{Y}$  is {\em full} if its monodromy group is the entire alternating or
 symmetric group on the degree.    There are two relatively simple obstructions to \eqref{maincover1} being
 full.  One is associated to $G$ having a non-trivial outer automorphism group and we deal with it
 by replacing $\AHur_h$ by a quotient variety $\AHur_h^*$ also covering $\AConf_\nu$.  The other
 is associated to $G$ having a non-trivial Schur multiplier and we deal with it by 
 a decomposition $\AHur_h^* = \coprod_\ell \AHur_{h,\ell}^*$.   Here $\ell$ runs over 
  an explicit quotient set of the Schur multiplier and each $\AHur_{h,\ell}^*$ is a union of connected components.

        The most important direction of the full-monodromy theorem is I $\Rightarrow$ II.  
        When $G$ is nonabelian and simple, this direction 
 is as follows.    \medskip
 
 \noindent {\em Fix a nonabelian simple group $G$ and 
 a list $C = (C_1,\dots,C_r)$ of conjugacy classes generating $G$. 
 Consider varying allowed $\nu$ and thus varying Hurwitz parameters $h = (G,C,\nu)$.  Then 
 as soon as $\min_i \nu_i$ is sufficiently large, 
  the covers $ \AHur^*_{h,\ell} \rightarrow \AConf_{\nu}$  are full and pairwise non-isomorphic.  }
\medskip

\noindent      The complete implication I $\Rightarrow$ II is similar, but 
$G$ is allowed to  be ``pseudosimple'', and therefore groups such
as $S_d$ are included.   There are considerable 
 complications arising from non-trivial abelianizations $G^{\rm ab}$,
 even in the case $|G^{\rm ab}| = 2$.  
 The extra generality is required for obtaining the natural converse
  II $\Rightarrow$ I.

 Our proof of I $\Rightarrow$ II in general starts from the Conway-Parker   theorem
about  connectivity of Hurwitz covers
   \cite{CP, FV,MM, EVW2}.   We deal with complications from
   nontrivial  $G^{\rm ab}$ in the framework of comparing two 
    Hochschild-Serre five-term exact sequences.  
   We upgrade connectivity to fullness by using a Goursat lemma
   adapted to our current situation and the explicit classification
   of $2$-transitive groups.  Our general approach has much in
   common with the proof of Theorem~7.4 in \cite{DT}, which is in
   a different context.  
   
   While there is a substantial literature on Hurwitz covers, our topic of asymptotic 
   fullness has not been systematically pursued before.   In related  directions
   there are the papers \cite{EEHS, MSV, Kluit}. We will indicate relations
   with some of this literature at various points in this paper.

\subsection{Specialization to number fields} \label{introNF}
 We say that a Hurwitz parameter $h = (G,C,\nu)$ is
{\em strongly rational} if all classes $C_i$ are rational. 
For strongly rational Hurwitz parameters,
\eqref{maincover1}  descends to a covering  
\begin{equation}
\label{maincover2}
\pi_h : \Hur_h \rightarrow \Conf_\nu
\end{equation}
 of $\Q$-varieties.    The full-monodromy theorem says 
 that every nonabelian finite simple group $G$ is part of infinitely many 
 Hurwitz parameters $h$ leading to full covers of $\Q$-varieties of 
 the modified form $\pi : \Hur^*_{h,\ell} \rightarrow \Conf_\nu$.
   
The group  $\Gal(\overline{\Q}/\Q)$ acts 
on  fibers $\pi^{-1}(u)$  over rational points $u \in \Conf_{\nu}(\Q)$.
The fullness of $\pi$, together with 
the Hilbert irreduciblity theorem, 
says that for {\em generic} $u$, the image of $\Gal(\overline{\Q}/\Q)$
contains the full alternating group on the fiber.  In this case one
gets a full field $K^*_{h,\ell,u}$ corresponding to the fiber.   

There is
a natural action of $\PGL_2(\Q)$ on $\Conf_\nu(\Q)$ and,
if $u,u'$ lie in the same $\PGL_2$-orbit, then $K^*_{h,\ell,u} \simeq K^*_{h,\ell, u'}$. 
Another application of the Hilbert irreducibility theorem
shows that, generically, different orbits give nonisomorphic fields.

Let $\cP$ be the set of prime divisors of $|G|$. 
The general theory of algebraic fundamental groups  
says that the action of $\Gal(\overline{\Q}/\Q)$ on 
$\pi^{-1}(u)$ is unramified away from $\primes$ 
so long as $u$ is a $\primes$-integral point.  
As $\nu$ varies, the number of $PGL_2$ equivalence classes of 
such specialization points can be arbitrarily large.  
 Thus,  {\em so long as even a weak 
 version of Hilbert irreducibility remains
 valid for such $\primes$-integral points}, we obtain 
 sufficiently many full fields for Conjecture~\ref{mc}.  
 As we explain in our last section, the
 available evidence is that there 
 is in fact a strong tendency for specializations
 to be full, and specializations from different
 orbits to be distinct.  
 
Of course, there are many other coverings of varieties
 that, like \eqref{maincover2}, have full geometric monodromy. However, 
 examples with the favorable ramification properties of the Hurwitz covers are rare.

 \subsection{Acknowledgements}
 We thank Simon Rubinstein-Salzedo and John Voight
  for helpful comments on earlier drafts of this paper. DR was supported by grant  
 \#209472  from the Simons foundation.  AV was supported by a grant from the NSF
 and from the Packard foundation. 
 
\section{Hurwitz covers}
     In this section we summarize the  theory of  Hurwitz covers,  taking the
     purely algebraic point of view necessary 
  for the application to Conjecture~\ref{mc}.   We consider Hurwitz parameters
  $h = (G,C,\nu)$, with $G$ assumed centerless to avoid technical complications.  
   The central focus is an associated cover
  $\pi_h : \Hur_h \rightarrow \Conf_\nu$ and related objects.
   A more detailed summary is in 
  \cite{RW} and a comprehensive reference is \cite{BR}.

\label{general}  

\subsection{Configuration spaces $\Conf_\nu$} 
\label{configspace} 
Let $\nu = (\nu_1,\dots,\nu_r)$ be a vector of positive
integers; we write $|\nu| = \sum \nu_i$.  For $k$ a field, let $\Conf_\nu(k)$ be the set of 
tuples $(D_1, \dots, D_r)$ of disjoint $k$-rational divisors on $\bbP^1_k$ with $D_i$ consisting of
$\nu_i$ distinct geometric points.    

Explicitly,   we may regard
\[
\Conf_\nu \subseteq \bbP^{\nu_1} \times \cdots \times \bbP^{\nu_r},
\]
where we regard $\bbP^{\nu_i}$ as the projectivized
space of binary homogeneous forms $q(x,y)$ of degree $\nu_i$,
and 
  $\Conf_{\nu}$ is then the open subvariety defined by nonvanishing
of the discriminant
$\disc(q_1 \cdots q_r)$.
The divisor $D_i$ associated to an $r$-tuple $(q_1,\dots,q_r)$ of such forms
is simply the zero locus of $q_i$.

\subsection{Standard Hurwitz varieties $\Hur_h$}  \label{Hurwitzscheme}

Let $k$ be an algebraically closed field of characteristic zero. 
Consider pairs $(\curve,\map)$ consisting of a proper smooth connected 
curve $\curve$ over $k$ together with a  Galois covering $\map: \curve \rightarrow \bbP^1$.

Such a pair has the following associated objects:
\begin{itemize}
\item An automorphism group $\Aut(\curve/\bbP^1)$ of size equal to 
the degree of $\map$, 
\item A branch locus $Z \subset \bbP^1(k)$ of degree $n= |Z|$;
 \item For every $t \in Z$, a local monodromy element $g_t \in \Aut(\curve/\bbP^1)$ defined up to conjugacy.  (To define this requires
 a compatible choice of roots of unity, i.e.\ an element of $\varprojlim_n \mu_n(k)$; we assume such a choice has been made). \end{itemize} 
  Consider triples
 $(\curve,\map,\iso)$ with $\iso: G \rightarrow \Aut(\curve/\bbP^1)$ 
 a given isomorphism.  We say that 
 such a triple has type $h$ if $\sum {\nu_i}=n$ and for each
 $i$ there are exactly $\nu_i$ elements  $t \in Z$ such that 
 $g_t \in C_i$.   The branch locus $Z$ then 
 defines an element of $\Conf_\nu(k)$ is 
 a natural way.

 The theory of Hurwitz varieties implies that
 there exists a $\overline{\Q}$-variety  $\Hur_{h}$, equipped with an {\'e}tale map
\begin{equation} \label{piHCdef}  \pi_h: \Hur_{h} \rightarrow \Conf_{\nu}, \end{equation} 
with the following property holding for all $k$:
 For any $u \in \Conf_{\nu}(k)$,  
 the fiber $\pi_h^{-1}(u)$ is, $\Aut(k/\Q(\mu_{\infty}))$-equivariantly, in bijection with the set of isomorphism classes
 of covers of $\PP^1$ of type $h$, with branch locus equal to $u$. 
 
  \subsection{Quotiented Hurwitz varieties $\Hur_h^*$} 
 If $(\curve, \map, \iso)$ is as above, 
 we can modify $\iso$ by an element $\alpha \in \Aut(G)$, to obtain a new triple $(\curve, \map, \iso \circ \alpha^{-1})$. 
 If $\alpha$ is inner, the resulting triple is actually isomorphic to $(\curve, \map, \alpha)$.
 As a result we obtain actions not of groups of automorphisms, but rather groups of outer automorphisms.

   Let $\Aut(G, C)$ be the subgroup consisting of those elements which fix every $C_i$. 
Then $\Out(G,C) = \Aut(G,C)/G$ acts naturally on $\Hur_h$ giving a quotient
 $$ \Hur^*_h = \Hur_h/ \Out(G, C),$$
 still lying over $\Conf_\nu$.  
 This quotient parameterizes pairs $(\curve, f)$ 
 equipped with an element $(D_1, \dots, D_r)$ of $\Conf_{\nu}(k)$
 so that the branch locus is precisely $\coprod D_i$, and there 
  exists an isomorphism
 $\iso: G \rightarrow \Aut(\curve/\PP^1)$ 
 so that the monodromy around each point of $D_i$ is of type $\iota(C_i)$.   
 Our main theorem focuses on $\Hur_h^*$ rather than $\Hur_h$. 
 
 \subsection{Descent to $\Q$} 
 \label{descent}
 The abelianized absolute Galois group 
 $\Gal(\overline{\Q}/\Q)^{\rm ab} = \hat{\Z}^{\times}$ 
 acts on the set of conjugacy classes in any finite group by 
 raising representing elements 
 to powers.   In particular, one can talk  about 
 ``rational'' classes, i.e.\ conjugacy classes fixed by this action.    As in the introduction, we say that 
  $h$ is {\em strongly rational}  if all $C_i$ are rational.  In this case, \eqref{piHCdef} 
  and its starred version $\pi_h^* : \Hur^*_h \rightarrow \Conf_\nu$ canonically descend to covers over $\Q$.  
  
  More generally, we say that $h$ is {\em rational} if conjugate classes appear
  with equal associated multiplicities.  In the main case when all the classes are different,
  this just means $\nu_i = \nu_j$ whenever $C_i$ and $C_j$
  lie in the same Galois orbit.   Rationality is a substantially weaker  
  condition than strong rationality.  For example, any finite group $G$ has rational $h$, 
  but only when $G^{\rm ab}$ is trivial or of exponent $2$ can $G$ have strongly
  rational $h$.
  
  For rational $h$, there is again canonical descent to $\Q$, although now the maps
  take the  form $\Hur_h \rightarrow \Hur_h^* \rightarrow \Conf_\nu^{\rho}$, with $\rho$ 
  indicating a suitable Galois twisting.         The subtlety of twisting is not seen
  at all in our main sections \S\ref{newtry2}-\S\ref{reverse}. 
 Our
  purpose in briefly discussing twisting here is to make clear that a large subset of 
  the covers considered in \S\ref{newtry2}-\S\ref{reverse} are useful in constructing
  fields for Conjecture~\ref{mc}.

  \section{Braid groups}
  \label{braid}  
       In this section we switch to a group-theoretic point of view, describing the monodromy 
  of Hurwitz covers $\pi_h : \AHur_h \rightarrow \AConf_\nu$  and $\pi^*_h : \AHur^*_h \rightarrow \AConf_\nu$ in
   terms of braid groups and their actions on explicit sets.   
   General references for braid groups and their monodromy actions
   include \cite[Chapter 3]{MM} and \cite[\S 2]{EEHS}. 
   
    Our main theorem concerns these monodromy 
   representations only, i.e.\ it is a theorem in pure topology. 
  As a notational device, used in the introduction and just now again,  
 we denote complex points by a different font as in $\AHur_h = \Hur_h(\C)$
 and $\AConf_{\nu} =\Conf_{\nu}(\C)$. 
 
 \subsection{Braid groups $\Br_\nu$}
  The Artin braid group on $n$ strands is defined by generators and relations:
\[
\Br_n = \left< \sigma_1,\dots,\sigma_{n-1} :
 \begin{array}{rclrcl} \sigma_i \sigma_j & = & \sigma_j \sigma_i,  & \mbox{if } |i-j| & > & 1 \\
 \sigma_i \sigma_j \sigma_i & = & \sigma_j \sigma_i \sigma_j, & \mbox{if } |i-j| & = & 1 \\
  \end{array} \right>.
\]
The rule $\sigma_i \mapsto (i, i+1)$ extends to a surjection $\Br_n \twoheadrightarrow S_n$.
For every subgroup of $S_n$, one gets a subgroup of $\Br_n$ by pullback.  
In particular, from the last component $\nu = (\nu_1,\dots,\nu_r)$ of a Hurwitz parameter one
gets a subgroup $S_\nu  := S_{\nu_1} \times \cdots \times S_{\nu_r}$.  We denote its
pullback by $\Br_\nu$.  The extreme $\Br_n$ above and the other extreme $\Br_{1^n}$ 
play particularly prominent roles in the literature, the latter often being called the
colored or pure braid group.

 \subsection{Fundamental groups} 
Let 
 $\star = (1,\dots,n) \in \AConf_{1^n}$. We will use it as a basepoint.
 We use the same notation $\star$ for its image in $\AConf_{\nu}$ for any $\nu$. 
  There is a standard surjection
$ \Br_{n} \twoheadrightarrow \pi_1(\AConf_{n}, \star)$
with kernel the smallest normal subgroup containing 
$\sigma_1 \cdots \sigma_{n-2} \sigma_{n-1}^2 \sigma_{n-2} \cdots \sigma_1$ \cite[Theorem~III.1.4]{MM}. 
 This map identifies $\sigma_i$
 with a small loop in  $\AConf_n$ that swaps the points $i$ and $i+1$.   
 Because of this very tight connection, the group $\pi_1(\AConf_n,\star)$  
 is often called the spherical braid group or the Hurwitz braid group.

     Similarly, we have surjections
    \begin{equation} \label{nuiso}  \Br_{\nu} \twoheadrightarrow \pi_1(\AConf_{\nu}, \star). \end{equation}
    Let $\cF_h$ and $\cF_h^*$ be the fibers of $\AHur_h$ and $\AHur^*_h$ over $\star$.   To completely translate into
    group theory, we need group-theoretical descriptions of these fibers as $\Br_\nu$-sets. 
     The remainder of this section accomplishes
    this task.

\subsection{Catch-all actions}  \label{catchall} We use the standard notational convention $g^h = h^{-1} g h$.  
 If $G$ is any group then $\Br_n$ acts on $G^n$ by means of a braiding rule,
  whereby $\sigma_i$ substitutes $g_i \rightarrow g_{i+1}$ and $g_{i+1} \rightarrow  g_i^{g_{i+1}}$:
     \begin{equation}
\label{bdef}
(\dots,g_{i-1},g_i,g_{i+1},g_{i+2},\dots)^{\sigma_i} = (\dots,g_{i-1},g_{i+1}, g_i^{g_{i+1}}, g_{i+2},\dots).
\end{equation}
Also $\Aut(G)$ acts on $G^n$ diagonally: 
\begin{equation}
\label{gdef}
(g_1,\dots,g_n)^\alpha = (g_1^\alpha, \dots, g_n^\alpha).  
\end{equation}
The braiding action and the diagonal action commute, so one has an action of the
product group $\Br_n \times \Aut(G)$ on $G^n$.

 \subsection{The $\Br_\nu$-sets $\mathcal{F}_h$ and $\cF_h^*$}  Next we replace $G^n$ by a smaller
 set appropriate to a given Hurwitz parameter $h$.  
  \label{Fhs} 
This smaller set is 
\begin{eqnarray}  \nonumber \mathcal{G}_h & = &  \{(g_1,\dots,g_n) \in G^n  :    g_1  \cdots g_n = 1, 
   \langle g_1, \dots,  g_n \rangle = G,  \\ &&  \label{Fhdef}  \qquad  \mbox{ first $\nu_1$ of the $g_i$'s lie in $C_1$, next $\nu_2$ lie in $C_2$, etc.}
\} . \end{eqnarray}
The subset $\cG_h$ is not preserved by all of $\Br_n \times \Aut(G)$, but it is preserved by 
$\Br_\nu \times \Aut(G,C)$.  The fibers then have the following group-theoretic description:
\begin{alignat}{2}
\label{id1} \cF_h &= \cG_h/\Inn(G) & & \simeq \mbox{(fiber of $\AHur_h \rightarrow \AConf_{\nu}$ above $\star$)},  \\
\label{id2} \cF_h^* &= \cG_h/\Aut(G,C) & &\simeq  \mbox{(fiber of $\AHur_h^* \rightarrow \AConf_{\nu}$ above $\star$)}.
\end{alignat}
Here in both cases the   isomorphisms $\simeq$ are isomorphisms of $\Br_{\nu}$-sets. 
Note that $\cF_h^* = \cF_h/\Out(G,C)$.   

\subsection{The asymptotic mass formula}   
\label{massformula} 
Character theory gives exact formulas for degrees, called mass formulas \cite{SeDar}.  We need 
only the asymptotic versions of these mass formulas, which are very simple:  
\begin{align}\label{ff8}
|\cF_h| & \sim \prod_{i=1}^r \frac{|C_i|^{\nu_i}}{|G'| |\Inn(G)|}, & 
|\cF^*_h| & \sim \prod_{i=1}^r \frac{|C_i|^{\nu_i}}{|G'| |\Aut(G,C)|}.
\end{align}
Here the meaning in each case is standard:   the left side over the right side  tends to $1$ for any
sequence of allowed $\nu$ with $\min_i \nu_i$ tending to $\infty$.    The structure of the products
on the right directly reflects the descriptions of the sets in \S\ref{Fhs}.

\section{Lifting invariants}
\label{lifting}
In this section we summarize the theory of lifting invariants which plays a key 
role in the study of connected components of Hurwitz spaces.    
  Group homology appears prominently and as a standing convention, we abbreviate   $H_i(\Gamma,\Z)$
by $H_i(\Gamma)$.

     In brief summary, the theory being reviewed goes as follows.  
     Let $h = (G, C,\nu)$ be a  Hurwitz \familial parameter.   
     The group $G$ determines its Schur multiplier $H_2(G)$. 
        In turn, $C$ determines
 a quotient group $H_2(G, C)$ of $H_2(G)$, and finally  $\nu$ determines a certain {torsor}
 $H_h = H_2(G, C,\nu)$ over $H_2(G, C)$.     The Conway-Parker theorem says
 that the natural map $\pi_0(\AHur_h) \rightarrow H_h$ is bijective whenever
 $\min_i \nu_i$ is sufficiently large.

 \subsection{The Schur multiplier $H_2(G)$} 
 \label{Schur}
 A stem extension of $G$ is a central extension ${G}^*$
such that the kernel of ${G}^* \rightarrow G$ is in the derived group of $G^*$.    A stem
extension of maximal order has kernel canonically isomorphic to the cohomology group $H_2(G)$.
This kernel is by definition the Schur multiplier.
A stem extension of maximal order is called a Schur cover.  
A given group can have non-isomorphic Schur covers, 
but this ambiguity  never poses problems for us here.

 \subsection{The reduced Schur multiplier $H_2(G, C)$}  \label{reducedSchur}
 
 If $x, y$ are commuting elements
of $G$, they canonically define an element $\langle x, y \rangle \in H_2(G)$:
the commutator of lifts of $x,y$ to a Schur cover.  This pairing is independent of the 
choice
of Schur cover.  In fact, a more intrinsic description is that $\langle x,y \rangle$ is the push forward
of the fundamental class of $H_2(\Z^2)$ under the map $\Z^2 \rightarrow G$
given by $(m,n) \mapsto x^m y^n$. 
 
 Fix a stem extension of maximal order $\tilde{G} \rightarrow G$.   For a conjugacy class
 $C_i$ and a list of conjugacy classes $C = (C_1,\dots,C_r)$ respectively, 
 define subgroups of the Schur multiplier:
 \begin{eqnarray}
\label{Ci} H_2(G)_{C_i} & = & \{ \langle g,z \rangle : g \in C_i \mbox{ and } z \in Z(g)\}, \\
\label{C} H_2(G)_C & = &  \sum H_2(G)_{C_i}.
\end{eqnarray}
  Here $Z(g)$ denotes the centralizer of $g$. 
The reduced Schur multiplier is then the corresponding quotient group.  
$H_2(G,C) = H_2(G)/H_2(G)_{C}$.

A choice of Schur cover $\tilde{G}$ determines a reduced Schur cover 
$\tilde{G}_C = \tilde{G}/H_2(G)_C$.  The corresponding short 
exact sequence 
  \[
  H_2(G, C) \hookrightarrow \tilde{G}_C \twoheadrightarrow G
  \]
 plays an essential role in our study.  
  
  In a degree $d$ central extension $\pi : G^* \rightarrow G$, the
  preimage of a conjugacy class $D$ consists of a certain number
  $s$ of conjugacy classes, all of size $(d/s) |D|$.  
    Always $s$ divides $d$.  If $s=d$ then $D$ is called {\em split.}  
    By construction, all the $C_i$ are split in $\tilde{G}_C$, and
    $\tilde{G}_C$ is a maximal extension with this property.  
    For more information on reduced Schur multipliers, see  \cite[\S 7, v1]{EVW2}.

\subsection{Torsors $H_2(G, C, \nu)$}   \label{subsec:torsors}

   For $i = 1$, \dots, $r$, let
 $H_2(G, C,i)$ be the set of conjugacy classes of $\tilde{G}_C$ that lie
 in the preimage of the class $C_i$.  If $\tilde{z}$ and $\tilde{g}$
 are lifts to $\tilde{G}_C$ of the identity $z=1$ and $g \in C_i$ respectively, then
 one can multiply $\tilde{z} \in H_2(G, C)$ and 
 $[\tilde{g}] \in H_2(G, C,i)$ to get 
$ [\tilde{z} \tilde{g}] \in 
  H_2(G, C,i)$.  This multiplication operator
  turns each $H_2(G, C,i)$
  into a torsor over $H_2(G, C)$.

One can multiply torsors over an abelian group: if $T_1$ and $T_2$ are torsors over an abelian group $Z$, 
then their product is $(T_1 \times T_2)/Z$ where all $(zt_1,z^{-1}t_2)$ have been identified. In our setting, one has a torsor 
\begin{equation} \label{torsordef}
 H_h := H_2(G, C,\nu) = \prod_i H_2(G, C,i)^{\nu_i}.
 \end{equation} 
 Note that $H_h$ 
 is naturally identified with the trivial torsor if all $\nu_i$ are multiples of the
 exponent of $H_2(G, C)$.  Namely the product $\prod a_i^{\nu_i}$ is independent of choices 
 $a_i \in H_2(G,C,i)$, and gives a distinguished element of $H_2(G, C, \nu)$. 
 In particular,  this distinguished element is fixed under $\Aut(G, C)$ (see  \S \ref{funcmap} for a more detailed discussion
 of functoriality).

\subsection{The lifting map}  \label{funcmapold}

Suppose given $(g_1, \dots, g_n) \in \cG_h$. 
 Lift each $g_{i}$ to an element of $\tilde{g}_{i} \in \tilde{G}_C$ arbitrarily,  
 subject to the unique condition that the product of the
 $\tilde{g}_{i}$ is the identity:
 $$\tilde{g}_1 \cdots \tilde{g}_n= 1 \in \tilde{G}_C.$$
 Then each $\tilde{g}_{i}$ determines
 an element $[\tilde{g}_{i}] \in H_2(G, C,i)$.
 Their product is  an element
 $\prod [\tilde{g}_i]  \in H_2(G, C,\nu),$ independent of choices.  
 This product is moreover  unchanged if we replaced $(g_1, \dots, g_n)$
 by another element in its $\Br_\nu$-orbit, or if we replace $(g_1, \dots, g_n)$ by
 a $G$-conjugate.  Thus, keeping in mind the identification $\pi_0(\AHur_h) = \cF_h/\Br_\nu$
 from \eqref{id1},
  we have defined a function
  \begin{equation}
  \label{lift}
 \inv_h :  \pi_0(\AHur_h)
   \longrightarrow H_h.
  \end{equation}
  We refer to $\inv_h$ as the lifting invariant.   It has been extensively studied by Fried and Serre, cf.\ \cite{BF,Se}. 
  When a set decomposes according to lifting invariants, we indicate this decomposition
  by subscripts.  Thus, e.g., $\cF_h = \coprod \cF_{h,\ell}$
  and $\cG_h = \coprod \cG_{h, \ell}$.
 
 The map \eqref{lift} is equivariant with respect to the natural actions of $\Out(G,C)$ and so we
 can pass to the quotient.  Writing 
  $H_h^*=H_h/\Out(G,C)$,
 we obtain
\begin{equation} \label{liftstar} \inv_h^* :  \pi_0(\AHur^*_h) \rightarrow H_h^*.\end{equation}
Again we notationally indicate lifting invariants by subscripts, so that 
e.g.\  $\mathcal{F}_{h,\ell}^* = \mathcal{F}_{h,\ell}/\Out(G,C)_{\ell}$, where  $\Out(G,C)_{\ell}$ is the stabilizer
 of $\ell$ inside $\Out(G,C)$.

Note that algebraic structure is typically lost in the process from passing from objects to their
corresponding starred objects.  Namely at the unstarred level, one has a group $H_2(G,C)$ and its many torsors
$H_h$.   At the starred level, 
$H^*_2(G,C)$ is typically no longer a group, the sets $H^*_h$ are no longer torsors, 
 and the cardinality of $H^*_h$ can depend on $\nu$.  Our main theorem makes direct reference
 only to $H_h^*$.  However in the proof we systematically lift from $H_h^*$ to $H_h$, to 
 make use of the richer algebraic properties.

We finally note for later use that there are asymptotic mass formulas for $\cF_{h,\ell}$
and $\cF_{h,\ell}^*$ that are very similar  to \eqref{ff8}. Indeed, they are derived simply
by applying \eqref{ff8} to $\tilde{G}_C$ together with liftings of the conjugacy classes $C_i$:
\begin{align}\label{refined-mass}
|\cF_{h,\ell}| & \sim  \frac{|\cF_h|}{|H_2(G,C)|}, & 
|\cF^*_{h,\ell}| & \sim \frac{|\cF_{h,\ell}|}{|\Out(G,C)_{\ell}|}.
\end{align}

 \subsection{Functoriality} \label{funcmap}

Suppose given a surjection $f: G \rightarrow H$
 of groups, together with conjugacy classes $C_i$ in $G$; set $D_i = f(C_i)$. This clearly
 induces a map $H_2(G, C) \rightarrow H_2(H,D)$. The functoriality of the torsors
 is less obvious,  because of the lack of uniqueness in a Schur cover. For this, we use a 
 more intrinsic presentation: 
 
 Amongst central extensions $\tilde{G} \rightarrow G$ equipped
 with a lifting $\tilde{C_i}$ of each $C_i$, there is a universal one $\tilde{G}^*$, unique up to unique isomorphism \cite[Theorem 7.5.1]{EVW2}.
 Now consider the central extension $G \times \Z^r \rightarrow G$, where we lift $C_i$ to $C_i \times e_i$, 
 where $e_i$ is the $i$th coordinate vector. This gives a canonical map $\alpha: \tilde{G}^* \rightarrow G \times \Z^r$, and we define
 $H_2(G, C, \nu)_{\univ}$ to be the preimage of $e \times \nu \in G \times \Z^r$. 
 
This is closely related to the previous definition. Note 
that if we fix lifts $C_i^* \subset \tilde{G}_C$ of each $C_i$, we get
an induced map $\beta: \tilde{G}^* \rightarrow \tilde{G}_C$ from the universal property.  This induces
a bijection $H_2(G, C, \nu)_{\univ}$ with $H_2(G,C)$; indeed, the canonical map
$$\alpha \times_G \beta: \tilde{G}^* \rightarrow \tilde{G}_C \times_G \left( G \times \Z^r \right) $$
is an isomorphism (again, \cite[Theorem 7.5.1]{EVW2}). 
So a choice of lifts $C_i^*$ give a distinguished element $c_{\nu} \in H_2(G, C, \nu)_{\univ}$ -- the preimage of the identity in $H_2(G,C)$.
 Moreover, 
if we replace $C_i^*$ by $z_i C_i^*$, where $z_i \in H_2(G, C)$, 
then the associated map $\tilde{G}^* \rightarrow \tilde{G}_C$ is multiplied by the composite map $\tilde{G}^* \rightarrow \Z^r \rightarrow \tilde{G}_C$
where the second map sends $e_i \in \Z^r$ to $z_i$. 
Thus, with this replacement, the identification $H_2(G, C, \nu)  \stackrel{\sim}{\rightarrow} H_2(G, C)$
has been multiplied by $z^{\nu_i}$; in other words, the distinguished element is replaced by $\prod z_i^{-\nu_i} c_{\nu}$.

This construction exhibits an identification of torsors
\begin{equation} \label{matchid} H_2(G, C, \nu)_{\univ} \simeq H_2(G, C, \nu)^{-1}, \end{equation}
where we write $T_1 \simeq T_2^{-1}$ for two $A$-torsors if there is an   identification of $T_1$ and $T_2$
transferring the $A$-action on $T_1$ to the inverse of the $A$-action on $T_2$. 

  In fact with respect to the identification \eqref{matchid},   our lifting invariant corresponds to the lifting invariant of \cite{EVW2}:
  In \cite{EVW2}, one takes $(g_1, \dots, g_r)$ and associates to it the lifting invariant $\Pi =\prod \widetilde{g_i} \in H_2(G, C, \nu)_{\univ}$, 
  where $\widetilde{g}$ is the lift to a universal central extensions equipped with lifting.  Fix $\tilde{G}_C$ and $C_i^*$
  and a morphism $\tilde{G}^* \rightarrow \tilde{G}_C$ as above.  Choose $z_i \in H_2(G, C)$ such that the image of $\Pi$ in $H_2(G, C)$ coincides with $\prod z_i^{\nu_i}$.
  Then $\prod \widetilde{g}_i $ is carried to $\prod z_i^{\nu_i}$ multiplied by the distinguished element of $H_2(G, C, \nu)_{\univ}$. On the
  other hand, the lifting invariant as we have defined it above equals 
  $[C_i^* z_i^{-1}] \in H_2(G, C, \nu)$, which equals $\prod z_i^{-\nu_i}$ times the corresponding element of $H_2(G, C, \nu)$.

Now -- returning to the surjection $G \rightarrow H$ --  take a universal extension $\tilde{H}^* \rightarrow H$ equipped with a lifting of the $D_i$, 
and consider $G \times_H \tilde{H}^* \rightarrow G$; it's a central extension
and it is equipped with a lifting of $C_i$, namely, $C_i \times_H D_i^*$.  There 
is thus a canonical map $\tilde{G}^* \rightarrow \tilde{H}^*$.   Taking fibers above $\nu \in \Z^r$ 
gives the desired map 
$$f_*: H_2(G, C, \nu)_{\univ} \rightarrow H_2(H, D, \nu)_{\univ},$$ 
and by inverting one obtains the desired map $H_2(G, C, \nu) \rightarrow H_2(H, D, \nu)$.  
   In particular, one easily verifies that if $H=G$ and $G \rightarrow H$ is an inner automorphism,
  the induced map on $H_2(G, C, \nu)$ is trivial. 
  
 Finally, suppose $\nu$ is chosen to be simultaneously divisible by the order of $H_2(G , C)$ and $H_2(H, D)$ (i.e.,
  each $\nu_i$ is so divisible.)   Then in fact the map 
  $H_2(G, C, \nu) \rightarrow H_2(H, D, \nu)$  respects the natural identifications of both
  sides with $H_2(G, C)$ and $H_2(H,D)$ (see after \eqref{torsordef}). In fact, one has  natural identifications
  $$H_2(G, C, \nu_1+\nu_2) \simeq H_2(G, C, \nu_1) \times H_2(G, C, \nu_2)/H_2(G,C),$$
  where the action of $z \in H_2(G,C)$ on the right is as $z : (t_1, t_2) \mapsto (t_1 z, z^{-1} t_2)$. 
These identifications are easily seen to be compatible with the map $H_2(G, C, \nu) \rightarrow H_2(H, D, \nu)$. 
Now choose $C_i^*$ and $D_i^*$ as above, giving rise to corresponding elements $c_{\nu} \in H_2(G, C, \nu)$
  and $d_{\nu} \in H_2(H, D, \nu)$. Write $f_* c_{\nu} = \gamma_{\nu} d_{\nu}$ for some $\gamma_{\nu} \in H_2(H, D)$;
  then our comments show that $\gamma_{\nu_1+\nu_2} = \gamma_{\nu_1} \gamma_{\nu_2}$, 
  and the claim follows: if $\nu$ is divisible by the order of $H_2(H, D)$, 
  then $\gamma_{\nu}$ will be trivial.

 \subsection{The Conway-Parker theorem}
     \label{conwayparker} 
  We will use a result due to Conway and Parker \cite{CP} in the important special case
  where $H_2(G,C)$ is trivial, and described in the paper of Fried and  V{\"o}lklein \cite{FV}.  See also
  \cite{EVW2, MM} for further information. 
   \begin{Lemma}  \label{cp} {\bf (Conway-Parker theorem)} Consider Hurwitz parameters $h = (G, C,\nu)$   
for $(G,C)$ fixed and $\nu$ varying.  Suppose that all the $C_i$ are distinct. 
  For sufficiently large $\min_i \nu_i$, the lifting invariant map $\inv_h : \pi_0(\AHur_h) \rightarrow H_h$ is bijective.  
\end{Lemma}
\noindent The Conway-Parker theorem plays a central role in this paper and a number of comments in several categories
are in order.  

First, the condition that $\min_i \nu_i$ is sufficiently large carries on passively to 
many of our later considerations.  We will repeat it explicitly several times but
also refer to it by the word {\em asymptotically.}     
 
Second, there are a number of equivalent statements.  The direct translation
of the bijectivity of $\pi_0(\AHur_h) \rightarrow H_h$  into group theory is that each fiber of 
 $\cF_h \rightarrow H_h$ is a single orbit of $\Br_\nu$.  Statements in the literature
 often compose the cover $\AHur_h \rightarrow \AConf_\nu$ with the 
 cover $\AConf_\nu \rightarrow \AConf_n$ and state the result in terms of actions
 of the full braid group $\Br_n$.   
 
 Third, quotienting   by $\Out(G,C)$ one gets a similar statement: 
 the resulting map $\inv_h^* : \pi_0(\AHur_h^*) \rightarrow H_h^*$ is asymptotically bijective.  
 This is the version that our full-monodromy theorem refines for certain $(G,C)$. 
  Note that a complication not present in Lemma~\ref{cp} itself appears at 
  this level:  the cardinality of $\cH^*_h$ can be dependent on $\nu$.  
 
\section{The full-monodromy theorem}
\label{newtry2}
     In this section, we state the full-monodromy theorem.  Involved in the statement is 
     a homological condition.  We clarify the nature of this condition by giving 
     instances when it holds and instances when it fails.

   \subsection{Preliminary definitions}  In this section, we define
   notions of {\em pseudosimple}, {\em unambiguous}, and {\em quasi-full}.  
   All three of these notions figure prominently in the statement of
   full-monodromy theorem.
   
       We say that a centerless finite group $G$ is {\em \pseudosimple} if 
       its derived group $\MZ$ is a power of a nonabelian simple group
  and any nontrivial quotient group of $G$ is abelian. 
 Thus, there is an extension
\begin{equation} \label{defext} \MZ \rightarrow G \rightarrow \Mab  \end{equation} 
     where $\MZ \simeq T^w$, with $T$ nonabelian simple, and the action of $\Mab$ on $T^w$
     is transitive on the $w$ simple factors.  [Our terminology is meant to be
     reminiscent of similar standard terms for groups closely related to a non-abelian
     simple group $T$:     {\em almost simple} groups are
     extensions $T.A$ contained in $\Aut(T)$ and {\em quasi-simple}
     groups are quotients $M.T$ of the Schur cover $\tilde{T}$.]
       
       We say that a conjugacy class $C_i$ in a group $G$ is {\em ambiguous} if
   the $G'$ action on $C_i$ by conjugation has more than one orbit.  If it has exactly 
   one orbit we say that $C_i$ is {\em unambiguous}.   These are standard notions
   and for many $G$ the division of classes into ambiguous and
   unambiguous can be read off from an Atlas page \cite{Atlas}.  
      
    Essentially repeating a definition from the introduction, we say that the action of a group $\Gamma$ on a set $X$ is {\em full}  
     if the image of $\Gamma$ in $\mathrm{Sym}(X)$ contains the alternating group $\Alt(X)$.
     Generalizing now, we say the action is {\em quasi-full} if the image contains $\Alt(X_1) \times  \cdots \times  \Alt(X_s)$, 
     where the $X_i$ are the orbits of $\Gamma$ on $X$.  Again we transfer the terminology to 
     a topological setting.  Thus a covering of a connected space $\mathsf{Y}$  is quasi-full 
     if for any $y \in \mathsf{Y}$, the monodromy action of $\pi_1(\mathsf{Y}, y)$  on
     $\mathsf{X}_y$ is quasi-full.   
   
     \subsection{Fiber powers of Hurwitz parameters}   This subsection describes how a Hurwitz parameter
     $h = (G,C,\nu)$ and a positive integer $k$ gives a triple $h^k = (G^{[k]},C^k,\nu)$.   
         Part of this
     notion, in the special case $k=2$, appears in the statement of the main theorem.  
     The general notion plays a central role in the proof. 
      
       In general, if $G$ is a finite group with abelianization $G^{\rm ab}$ we can consider its
      $k$-fold fiber power
      \[
      G^{[k]} = G \times_{G^{\rm ab}}  \cdots \times_{G^{\rm ab}} G.
      \]
      Note that even when $G = T^w.G^{\rm ab}$ is pseudosimple, the fiber powers $G^{[k]} = T^{wk}.G^{\rm ab}$ 
      for $k \geq 2$ are not, because $G^{\rm ab}$ does not act transitively on the factors.  
      
     If $C_i$ is a conjugacy class in a group $G$, we can consider its  Cartesian powers 
     $C_i^k \subseteq  G^{[k]}$.    In general, $C_i^k$ is only a union of conjugacy classes.
     However if $C_i$ unambiguous then $C_i^k$ is a single class.  
     
     If $C = (C_1,\dots,C_r)$ is a list of conjugacy classes, we can consider the 
     corresponding list $(C^k_1,\dots,C_r^k)$.   Generation of $G$ by the $C_i$ does 
     not imply generation of $G^{[k]}$ by the $C_i^k$.  
          However if $G$ is pseudosimple then this implication does hold.   Thus
     if $G$ is pseudosimple and $C$ consists only of unambiguous classes, 
     the triple $h^k$ is a Hurwitz parameter.    
    
    Suppose, then, that $G$ is pseudosimple and $C$ consists of unambigous classes.  
The natural map (\S \ref{funcmap}) 
     $$H_2(G^{[k]}, C^k,\nu) \rightarrow H_2(G, C, \nu)^k$$
      is surjective.       This surjectivity can be seen by interpreting both sides in terms of connected components via the Conway-Parker theorem.
     Alternately, it follows because the map is equivariant with respect
     to the natural map $H_2(G^{[k]}, C^k) \rightarrow H_2(G, C)^k$, 
     which is surjective by homological algebra, as we explain after 
\eqref{ei2}.

     \subsection{Statement}  With our various definitions in place, we can state the main result of this paper.  
     
      \begin{Theorem} {\bf (The full-monodromy theorem)} \label{mt}
Let $G$ be a finite centerless nonabelian group, let $C = (C_1,\dots,C_r)$ a list of  distinct nonidentity conjugacy classes generating
$G$,
and consider Hurwitz parameters $h = (G,C,\nu)$ for varying allowed $\nu \in \Z_{\geq 1}^r$.    
Then the following are equivalent:
\begin{description}
\item[I]
\begin{itemize}
\item[1.] $G$ is pseudosimple,
\item[2.] The classes $C_i$ are all unambiguous, and   
\item[3.] $|H_2(G^{[2]},C^2)| = |H_2(G,C)|^2$.
\end{itemize}
\item[II] The covers $\AHur^*_{h} \rightarrow \AConf_\nu$ are quasi-full whenever $\min_i \nu_i$ is 
sufficiently large.  
\end{description}
\end{Theorem}
\noindent

Note that Statement II can equivalently be 
presented in terms of fullness:  for
$\min_i \nu_i$ sufficiently large, the covers $\AHur_{h,\ell}^* \rightarrow \AConf_\nu$ 
are full and pairwise non-isomorphic as $\ell$ ranges over $H_h^*$.   
 Note also that a pseudosimple group 
$G$ is simple if and only if $G^{\rm ab}$ is trivial.   In this case, Conditions I.2 and I.3 are trivially satisfied and the direction I 
$\Rightarrow$ II becomes the statement highlighted in \S\ref{introFM}.

     For the  more important direction I $\Rightarrow$ II, the condition that $\min_i \nu_i$ is
 sufficiently large is simply inherited from the Conway-Parker theorem. 
  Calculations suggest that quasi-fullness tends to be obtained as soon 
 as it is allowed by the mass formula.
 We are not pursuing  the important question of effectivity here,  but we note that effective
 statements of fullness are obtained for certain classical Hurwitz parameters
 in \cite{Kluit}.

     Given $(G,C)$, whether or not Conditions~1 and 2 hold is immediately determinable
in practice.   Evaluating Condition~3 is harder in general, and the next two subsections
are devoted to giving an easily checkable reformulation applicable in many cases  (Proposition~\ref{homprop})
and showing  (Corollary~\ref{homcor}) that it sometimes fails.

\subsection{The homological condition for $G$ of split-cyclic type}  
We say that a pseudosimple group $G$ has {\em split} type
if the canonical surjection 
$\pi : G \rightarrow G^{\rm ab}$ has a section $s : G^{\rm ab} \rightarrow G$.  
 This {\em a priori} strong
condition is actually commonly satisfied. 
Similarly, we say that a pseudosimple group has {\em cyclic} 
type if $G^{\rm ab}$ is cyclic.    Again this strong-seeming condition is commonly satisfied, as indeed
for a simple group $T$ all of $\Out(T)$ is often cyclic \cite{Atlas}.  
 When both of these conditions
are satisfied, we say that $G$ is of {\em split-cyclic} type.  

For $G$ of split-cyclic type, the next proposition says that Condition~3 of 
Theorem~\ref{mt} is equivalent to an apparent strengthening $\hat{3}$.  
Moreover these two conditions are both equivalent to 
a more explicit condition E which makes
no reference to either fiber powers or powers.  For E, we modify
the notions defined in  \S\ref{reducedSchur} as follows:
\begin{eqnarray*}
H'_2(G)_{C_i} & = & \{ \langle g,z \rangle : g \in C_i \mbox{ and } z \in Z(g) \cap G'\}, \\
H'_2(G)_C & = & \sum H'_2(G)_{C_i}.
\end{eqnarray*}
These are straightforward variants, as indeed if one removes every $'$ one
recovers definitions \eqref{Ci}, \eqref{C} of the previous notions.  

\begin{Proposition} 
\label{homprop} Let $G$ be a pseudosimple group of split-cyclic type and
let $C = (C_1,\dots,C_r)$ be a list of distinct unambiguous conjugacy classes.
Then the following are equivalent:
\begin{description}
\item[$3$]  $|H_2(G^{[2]},C^2)| = |H_2(G,C)|^2$,
\item[$\hat{3}$] $|H_2(G^{[k]},C^k)| = |H_2(G,C)|^k$ for all positive integers $k$,
\item[$\mathrm{E}$] $H_2(G)_C = H'_2(G)_C$.
\end{description}
Moreover if $|G^{\rm ab}|$ is relatively prime to $|H_2(G)|$ then all
three conditions hold.  
\end{Proposition}

\proof   
 All three conditions involve the list $C$ of conjugacy classes.  We begin however with
 considerations involving $G$ only.  
 The $k$ different coordinate projections $G^{[k]} \rightarrow G$ together
 induce a map $f_k : H_2(G^{[k]} ) \rightarrow H_2(G )^k$.  
 We first show that the assumption that $G$ has split-cyclic type
 implies all the $f_k$ are isomorphisms.   We present this deduction in some detail
 because we will return to parts of it in \S\ref{Endgeneral}.
 
The map $f_k$ is part of a morphism of five-term exact sequences (see \cite[Theorem 5.2]{EckmannStambach}, noting that $H_1(G') = 0$)
{\renewcommand{\arraycolsep}{1pt}
\begin{equation}  \label{ei2}
\begin{array}{ccccccccc} H_3(\M^{[k]}  ) & \stackrel{\pi_3^{[k]}}{ \rightarrow} & H_3(\Mab  )  & \stackrel{\delta^{[k]}}{\rightarrow}  &H_2(\MZ^k  )_{\Mab} &  \stackrel{i^{[k]}_2}{\rightarrow}   &
H_2(\M^{[k]}  ) & \stackrel{\pi^{[k]}_2}{\twoheadrightarrow} & H_2(\Mab  )  \\
 \downarrow && \downarrow{\Delta_3} && \downarrow{\simeq} && \downarrow{f_k} && \downarrow{\Delta_2} \\
 H_3(\M  )^k & \stackrel{\pi_3^k}{ \rightarrow} & H_3(\Mab  )^k &  \stackrel{\delta^k}{\rightarrow} & H_2(\MZ  )^k_{\Mab} &
 \stackrel{i^k_2}{\rightarrow}  &
H_2(\M  )^k & \stackrel{\pi^k_2}{\twoheadrightarrow} &  H_2(\Mab  )^k. 
\end{array}
\end{equation} 
}

\noindent   Each five-term sequence arises from the Hochschild-Serre spectral sequence associated
to an exact sequence of groups.  The top sequence comes from the $k^{\rm th}$ 
fiber power of $G' \stackrel{i}{\hookrightarrow} G \stackrel{\pi}{\twoheadrightarrow} G^{\rm ab}$,
while the bottom sequence comes from the $k^{\rm th}$ ordinary Cartesian power.

We note that \eqref{ei2} actually shows that $H_2(G^{[k]}, C^k) \rightarrow H_2(G,C)^k$ is surjective
whenever $G$ is pseudosimple and $C$ consists of unambiguous classes. 
The point is that $H_2(G)_C$ surjects onto $H_2(G^{\ab})$. That is 
because $H_2(G^{\ab})$ is generated by symbols $\langle \alpha, \beta \rangle$.
But such a symbol belongs to the image of $H_2(G)_C$,
since the $[C_i]$ generate $G^{\ab}$ and, for any $g \in C_i$,
the centralizer $Z(g)$ surjects to $G^{\ab}$ because $C_i$ is unambiguous.

The  assumption that $\pi : G \rightarrow G^{\rm ab}$ has a splitting $s$ drastically simplifies \eqref{ei2}.    
From $\pi \circ s = \mbox{Id}_{G^{\rm ab}}$ 
one gets that  $\pi^{[k]}_3 \circ s^{[k]}_3$ and $\pi^k_3 \circ s^k_3$ are the identity on 
$H_3(\Mab  )$ and $H_3(\Mab  )^k$ respectively.  Thus $\pi^{[k]}_3$ and $\pi_3^{k}$ are both
surjective and so the boundary maps $\delta^{[k]}$ and $\delta^k$ are both $0$.    Thus the part of
\eqref{ei2} relevant for us becomes
\begin{eqnarray}
 \label{3term}
 \begin{array}{ccccc} 
 H_2(\MZ^k  )_{\Mab} & {\hookrightarrow}  &
H_2(\M^{[k]}  ) & \twoheadrightarrow & H_2(\Mab  ) \\
\downarrow{\simeq} && \downarrow{f_k} && \downarrow{\Delta_2} \\
H_2(\MZ  )_{\Mab}^k & {\hookrightarrow}  &
H_2(\M  )^k&  \twoheadrightarrow & H_2(\Mab  )^k.
\end{array}
\end{eqnarray}
We have suppressed some notation, since we have no further use for it.  

The assumption that $G^{\rm ab}$ is cyclic is equivalent to the assumption 
that $H_2(G^{\rm ab} )$ is the zero group.  Thus exactly in this situation
one gets the independent simplification of \eqref{ei2} 
where the last column becomes the zero map between zero groups.
Applied to \eqref{3term} it says that $f_k : H_2(G^{[k]} ) \rightarrow H_2(G )^k$
is an isomorphism.  We henceforth use $f_k$ to identify $H_2(G^{[k]} )$ with  $H_2(G )^k$.

We now bring in the list $C$ of conjugacy classes.  
We have a morphism of short exact sequences:
\begin{equation}
  \begin{array}{ccccc}
   H_2(G^{[k]} )_C & \hookrightarrow & H_2(G^{[k]} ) & \twoheadrightarrow & H_2(G^{[k]},C^k) \\
   \rotatebox{270}{$\!\!\!\! \subseteq$} & & \rotatebox{270}{$\!\!\!\! =$}   & &  \rotatebox{270}{$\!\!\!\! \twoheadrightarrow$} \\
    H_2(G )_C^k & \hookrightarrow & H_2(G )^k & \twoheadrightarrow & H_2(G,C)^k.
  \end{array}
  \end{equation} 
  Since the map in the right column is surjective, Conditions $3$ and $\hat{3}$ become 
  that it is an isomorphism for $k=2$ and all $k$ respectively.  So they are equivalent
  to the inclusion in the left column being equality, again
  for $k=2$ and all $k$ respectively.  We work henceforth with these
  versions of Conditions $3$ and $\hat{3}$.  
 
Trivially 
\begin{equation}
\label{ABC1}  |H_2(G)_C | = \left|  H_2'(G)_C \right| \cdot \left| H_2(G)_C/H_2'(G)_C \right|.
\end{equation}
But also the image of $H_2(G^{[k]})_{C^k}$ in $\left( H_2(G  )/H_2'(G)_C \right)^k$
is exactly the diagonal image of $H_2(G)_C/H_2'(G)_C$.
To see this, note that $H_2(G^{[k]})_{C^k}$  is generated by
$$(\langle g, z_1 \rangle, \dots, \langle g, z_k \rangle) $$
where $g \in \bigcup C_i$, each $z_i \in Z(g)$, and $z_1 \equiv \cdots \equiv z_k$ 
modulo $G'$.  In particular, it  certainly contains the diagonal image of $H_2(G)_C$. 
   On the other hand,
the images of $\langle g, z_i \rangle$ inside $H_2(G)_C/H_2'(G)_C$
are equal to each other, since $\langle g, z_i z_j^{-1} \rangle \in H_2'(G)_C$. 

Moreover, $H_2'(G)_C^k \subseteq H_2(G^{[k]})_{C^k}$. 
This inclusion holds because, for any $g \in C_i$ and $z \in Z(g) \cap G'$, we have
$$ (\langle g, z \rangle, 0, 0, \dots) \in H_2(G^{[k]})_{C^k},$$
since we can regard the left-hand side as
$(\langle g,z \rangle, \langle g, e \rangle, \langle g, e \rangle, \dots )$. 
Similarly for any other ``coordinate.'' 
Therefore,
\begin{equation}
\label{ABC2}  |H_2(G^{[k]})_{C^k} |  = \left| H_2'(G)_C \right|^k  \cdot \left|H_2(G)_C/ H_2'(G)_C \right|.
\end{equation}
Dividing the $k^{\rm th}$ power of \eqref{ABC1} by \eqref{ABC2}, one gets
\begin{equation}
\label{ABC}
 \frac{|H_2(G)_C|^k}{|H_2(G^{[k]})_{C^k}|}  = |H_2(G)_C/H_2'(G)_C|^{k-1}.
\end{equation}  
Condition $3$ says the left side is $1$ for $k=2$.  Condition $\hat{3}$ says the left side is $1$ for all $k$.  
Equation~\eqref{ABC} says that each of these is equivalent to $H_2(G)_C = H_2'(G)_C$, which
is exactly Condition E.

For the final statement,  $|H_2(G)_{C_i}/H_2'(G)_{C_i}|$ clearly divides
$|H_2(G)|$.  It also divides $|G^{\rm ab}|$, because
$Z(g)/(Z(g) \cap G')$ surjects onto $H_2(G)_{C_i}/H_2'(G)_{C_i}$
via $z \in Z(g) \mapsto \langle g, z \rangle$, for any fixed $g \in C_i$.
So, if $|H_2(G)|$ and $|G^{\rm ab}|$ are relatively prime then
always $H_2(G)_{C_i} = H_2'(G)_{C_i}$ and so Condition E holds.  
 \qed

\subsection{The homological condition for $G$ of split-$p$-$p$ type}
\label{splitpp}
For $p$ a prime, we say that a pseudosimple group $G$ has {\em split}-$p$-$p$ type if 
$G \rightarrow G^{\rm ab}$ is split and 
\[
|G^{\rm ab}| = |H_2(G )| = p.
\]
Even this seemingly  very special case is common.  For example, taking $p=2$,  it
includes 
\begin{itemize}
\item all six extensions $T.A$ of sporadic groups $T$ with $A$ and $H_2(T.A )$ 
non-trivial, 
\item all
$S_d$ with $d \geq 5$, and
\item  all $PGL_2(q)$ for  odd $q \geq 5$.   
\end{itemize}
To illustrate the tractability of Condition~E of Proposition~\ref{homprop}, we work it out explicitly 
for groups $G$ of split-$p$-$p$ type.   Explicating Condition~E  for the full split-cyclic case would be similar
but combinatorially  more complicated.  
 
     For  $G$ of split-$p$-$p$ type, we divides its unambiguous classes up into
 three types.  Let $\tilde{G}$ be a Schur cover of $G$.  An unambiguous class $C$ is {\em split} if its preimage
$\tilde{C}$ consists of $p$ conjugacy classes in $\tilde{G}$.  
It is {\em mixed} if 
 $\tilde{C}$ is 
 $p$ different $\tilde{G}'$ conjugacy classes but just one 
  $\tilde{G}$ class.  Otherwise a class $C$ is {\em inert}.   Mixed classes are necessarily in the derived group, but 
split and inert classes can lie above any element of $G^{\rm ab}$. 

\begin{Corollary}
\label{homcor} Let $G$ be a pseudosimple group of split-$p$-$p$ type
 and let $C = (C_1,\dots,C_r)$ 
be a list of unambiguous classes.  
Then Condition E fails exactly when the are no inert classes and at least one
mixed class among the $C_i$.  
\end{Corollary}

\proof  We are considering subgroups of the $p$-element Schur multiplier $H_2(G )$.  The subgroups have the following  
form
\[
\begin{array}{c|ccc}
C_i & \mbox{Split} & \mbox{Mixed} & \mbox{Inert} \\
\hline
H_2'(G)_{C_i} & 0 & 0 & H_2(G) \\
H_2(G)_{C_i} & 0 & H_2(G) & H_2(G)
\end{array}
\]
Thus $H_2'(G)_C = \sum_i H_2'(G)_{C_i}$ is a proper subgroup of $H_2(G)_C = \sum_i H_2(G)_{C_i}$ 
exactly under the conditions stated in the corollary. \qed 

    For a group $T.p$, the types of classes can be determined from an Atlas-style
character table, including its lifting row and fusion column.  For example, for the six sporadic
$T$ mentioned above, the mixed classes in $T.2$ are exactly as follows:  
\[
{\renewcommand{\arraycolsep}{4pt}
\begin{array}{c|c|c|c|c|c}
\mbox{Mathieu}_{12} & \mbox{Mathieu}_{22} & \mbox{Hall-Janko} & \mbox{Higman-Sims} & \mbox{Suzuki}& \mbox{Fischer}_{22}  \\
\hline
10A & 8A & 8A & 4A, 6A, 12A& 12D, 12E, 24A & (15 \mbox{ classes})
 \\ 
\end{array}
}
\]
In the sequences $S_d$ and $PGL_2(q)$, the patterns evident from character tables in the first few instances
can be proved to hold in general.   Namely for $S_d$, conjugacy classes are indexed by partitions of $d$.  
The type of a class $C_\lambda$ 
can be read off from two features of the indexing partition $\lambda$, the number $e$ of even parts and
whether or not all parts are distinct:
\[
\begin{array}{c|ccc}
  &  e=0&  e \in \{2,4,6,\dots\} & e \in \{1,3,5,\dots\}  \\
\hline
\mbox{All distinct} & \mbox{Ambiguous} & \mbox{Mixed} & \mbox{Split}  \\
\mbox{Not all distinct} & \mbox{Split} & \mbox{Inert} & \mbox{Inert}
\end{array}
\]
Thus $S_5$ has no mixed classes while $C_{42}$ and $C_{421}$ are the unique
mixed classes of $S_6$ and $S_7$ respectively.    For $PGL_2(q)$, the
division is even easier:  the two classes of order the prime dividing $q$ are ambiguous, the
two classes of order $2$ are inert, and all other classes are split.  Thus for
$PGL_2(q)$, the homological condition always holds.

\section{Proof of $\mbox{I} \Rightarrow \mbox{II}$}  
\label{direct}
    In this section we prove the implication $\mbox{I} \Rightarrow \mbox{II}$ of Theorem~\ref{mt}.  
 Thus we consider Hurwitz parameters $h = (G,C,\nu)$ for fixed $(G,C)$ satisfying Conditions~1-3 
 and varying $\nu$.  We then prove that the action of $\Br_\nu$ on $\cF^*_h$ is quasi-full
 whenever $\min_i \nu_i$ is sufficiently large.  The implication $\mbox{I} \Rightarrow \mbox{II}$
 is the part of Theorem~\ref{mt} which provides theoretical support for Conjecture~\ref{mc}.  
 
\subsection{A Goursat Lemma}
 
The classical Goursat lemma classifies certain subgroups of powers of a simple group.
We state and prove a generalized version here. As usual, if 
one has groups $G_1$, $G_2$ endowed with homomorphisms $\pi_1$, $\pi_2$ to a third 
group $Q$, we say that $G_1$ and $G_2$ are isomorphic {\em over $Q$} if
there is an isomorphism $i : G_1 \rightarrow G_2$ satisfying $\pi_2 i = \pi_1$.  
\begin{Lemma}   {\bf (Generalized Goursat lemma) }
\label{gs} 
Suppose that $\M$ is \pseudosimple,     and $H \subseteq  G^{[k]}$ is 
a ``Goursat subgroup'' in the sense that it surjects onto each coordinate factor.
Then 
\begin{itemize}
\item[1:]    $H$ is itself isomorphic over $G^{\ab}$ to $G^{[w]}$ for some $w \leq k$. 
\item[2:]    There is a surjection  $f: [1,k]  \rightarrow [1,w]$
and  automorphisms $\varphi_1$, \dots, $\varphi_k$  of $G$ over $G^{\ab}$ 
such that $H$ is the image of
$G^{[w]}$ under
$$ (g_1, \dots, g_w) \mapsto (\varphi_{1}(g_{f(1)}), \dots,  \varphi_k(g_{f(k)}) ).$$
\end{itemize}
  \end{Lemma}  
 \proof 
  
 We first prove Statement 1 by induction, the base case $k=1$ being trivial.  Note that the projection  $\bar{H} = \pi_2(H)$ of $H$ to the second factor in
 $$G^{[k]} = G \times_{G^{\ab}} G^{[k-1]}$$
 is also a Goursat subgroup. By induction, it is $G^{\ab}$-isomorphic
 to $G^{[v]}$ for suitable $v$.  The kernel $K = \mathrm{ker}(\pi_2)$ of  the projection $H \rightarrow \bar{H}$
  maps, under the first projection $\pi_1$, to a
 subgroup  $\bar{K} \subseteq \MZ$ that is invariant under conjugation by $\M$. 
 In particular, either $\bar{K}$ is trivial, and we're done by induction, 
 or $\bar{K} = \MZ$. In the latter case,  we will show
 that $H = G \times_{G^{\ab}} \bar{H}$: 
Take any element 
 $ (m^*, \mu) \in G \times_{G^{\ab}}  \bar{H}$. By assumption
 there exists $m$ in $G$ such that
 $ (m, \mu) \in H$; but then $m$ and $m^*$ have the same projection to $G^{\ab}$,
 and so 
 $$ (m^*, \mu) = (m^* m^{-1}, 1) \cdot (m , \mu)$$
 lies in $H$ also.   This concludes the proof of the first assertion: $H$ is isomorphic to $G^{[w]}$ over $G^{\ab}$
 for some $w$.

 Now we deduce Statement 2 from Statement 1.   
  Let $\Theta = G^{[w]} \rightarrow H$ be any isomorphism
  and write $\Theta(g) = (\theta_1(g),\dots,\theta_k(g))$.  
  We need to show that for each $i$, one can express
  $\theta_i(g)$ in the form $\varphi_i(g_{f(i)})$ as 
  in Statement $2$.   In other words, letting $\pi_j : G^{[w]} \rightarrow G$
  be the $j^{\rm th}$ projection, we need
  to show that any surjective morphism 
  $\theta :  G^{[w]} \rightarrow G$ over $G^{\ab}$ 
  factors as $\varphi \pi_j$ for some $j \in \{1,\dots,w\}$ and some 
  automorphism $\varphi : G \rightarrow G$ over $G^{\rm ab}$.

So let $\theta: G^{[w]} \rightarrow G$ be any surjective morphism 
over $G^{\ab}$.    Its kernel $K$ is a normal subgroup of $(G')^w$, invariant under $G^{[w]}$, 
 and with index $|G'|$.      Now, via $G' \simeq T^u$ for some nonabelian simple group $T$, 
  the normal subgroups of $(G')^w \simeq T^{uw}$
 are of the form $T_I = \prod_{(i,j) \in I} T_{(i,j)}$, where $I$ is a subset of 
 $P = \{1, \dots, u\} \times \{1,\dots, w\}$.   The normal subgroups which 
 are invariant under $G^w$ are those for which 
 the indexing set $I$ is invariant
 under the natural action of $G^{\rm ab}$.  
 The orbits of $G^{\rm ab}$ on $P$ are the sets $P_j =  \{1,\dots,u\} \times \{j\}$.  
 So the kernel $K$ of $\theta$ necessarily has the form $T_{P - P_j}$.
 Thus $K$ is also the kernel of the coordinate projection
  $\pi_j$.   The unique bijection $\varphi : G \rightarrow G$ 
  satisfying  $\theta = \varphi \pi_j$ is then an automorphism of 
  $G$ over $G^{\rm ab}$.  
  \qed

\subsection{Identifying braid orbits}  \label{core} 

\label{proof1}  For $F$ a set and $k$ a positive integer we let 
\[
F^{\underline{k}} = \{(x_1,\dots,x_k) : \mbox{all $x_i$ are different}\}.
\]
If $F$ has cardinality $N$ then $F^{\underline{k}}$ has cardinality 
$N^{\underline{k}} := N(N-1) \cdots (N-k+1)$.  
In this subsection we assume Conditions 1 and 2 and
identify the quotient set $(\cF_h^*)^{\underline{k}}/\Br_\nu$ 
asymptotically.    
  
Begin with $x_1, \dots, x_k \in \cF^*_h$. 
Choose a set of representatives $\tupleg_1, \dots, \tupleg_k \in \cG_h$.   Writing each $\tupleg_i$ as a column 
vector, we get a matrix
 \begin{equation}
 \label{gmatrix} 
( \tupleg_1, \dots, \tupleg_k ) =  \left( \begin{array}{cccc}  g_{11} & g_{21} & \dots & g_{k1} \\
 g_{12} & g_{22} & \cdots & g_{k2} \\
 g_{13} & g_{23} & \cdots & g_{k3} \\
 \vdots & \vdots & \vdots & \vdots \\
 g_{1n} & g_{2n} & \cdots & g_{kn}
 \end{array} \right).
 \end{equation}
 So, simply recalling our context:
 \begin{itemize}
 \item   All the $g_{ij}$ in  
 in a given row are in the same
 conjugacy class of $\M$. 
  \item These conjugacy classes are
 $\overbrace{C_1,\dots,C_1}^{\nu_1}; \dots; \overbrace{ C_r, \dots, C_r}^{\nu_r}$
 as one goes down the rows, so that a given row is 
 in some $C_i^k$.
\item  Each column in
 its given order multiplies to $1$.
 \item  Each column generates all of $\M$.  
 \end{itemize}
 
All entries in a given row certainly have the same projection to $G^{\rm ab}$ and so each
row defines an element of $G^{[k]}$.  Consider now the subgroup $H$ of $G^{[k]}$
 generated by the rows of this matrix.   We are going to show that
\begin{equation} \label{disjointgeneration} H=G^{[k]} \iff 
\mbox{(all $x_i$ are different)}.
 \end{equation}

First of all, note that the condition that $H=G^{[k]}$
  is independent of the choice of lifting from
$\cF_h^*$ to $\cG_h$. For example, if we modify the
$\tupleg_1$, the first column of \eqref{gmatrix},
by an element $\alpha \in \Aut(G,C)$, then the subgroup generated by the rows
simply changes by the automorphism $(\alpha, 1, 1,1 \dots, 1)$ of $G^{[k]}$. 
Note that $\alpha$ is automatically an isomorphism of $G$ over $G^{\ab}$
because it preserves each $C_i$ and they generate $G^{\ab}$.

Now direction $\implies$ of \eqref{disjointgeneration} is easy: if
$x_i = x_j$ for some $i \neq j$ then we could lift so that $\tupleg_i = \tupleg_j$, 
and then certainly $H \subsetneq G^{[k]}$. 

Now suppose that $x_i \neq x_j$
 for all $i \neq j$; we'll show that $H=G^{[k]}$.
 Since
 each column generates $\M$, the subgroup  $H$  is a
 Goursat subgroup of $G^{[k]}$.   Accordingly we may apply
  Lemma~\ref{gs}, and see that $H$
  can be constructed from a surjective function $f: [1,k] \rightarrow [1,w]$
  together with   a system of isomorphisms $\varphi_j: G \rightarrow G$
  over $G^{\ab}$, for $1 \leq j \leq k$. 
  In particular, 
 we may find $(y_1, \dots, y_w) \in G^{[w]}$
 which maps
 to the first row 
 $(g_{11}, g_{21}, \dots, g_{k1})$, so that
 $$ \varphi_{j} (y_{f(j)}) = g_{j1}, \ \ 1 \leq j \leq k.$$
 In particular, whenever $f(j) = f(j')$, 
 the map
 $$\varphi_{j'} \varphi_j^{-1}$$
carries $g_{j1}$ to $g_{j'1}$ and so  preserves $C_1$.   By similar reasoning, applied to the second row, third row and so on,  
this map  preserves
 {\em every} conjugacy class, so
 $$ \varphi_{j'} \varphi_{j}^{-1} \in \Aut(G, C)$$
 whenever $f(j) = f(j')$.  But  $\varphi_{j'} \varphi_{j}^{-1}$
 carries $g_{ji}$ to $g_{j' i}$; 
 that means that actually $x_j = x_j'$,
 and so $j=j'$. In other words, $f$ is injective, 
 and so $H \simeq G^{[k]}$ as desired.

Each matrix \eqref{gmatrix} with $H$ all of $G^{[k]}$ defines an element of $\cG_{h^k}$.  
Now, the group $\Aut(G,C)^k$ acts on $G^{[k]}$; its image in the outer automorphism group
will be called $\Out(G,C)^{[k]}$. This latter group maps onto $\Out(G,C)^k$, with kernel
isomorphic to $(G^{\ab})^{k-1}$.   Our considerations have given a bijective map
\begin{equation} 
\cF_{h^k}/\Out(G,C)^{[k]} \stackrel{\sim}{\longrightarrow} \cF_h^{*\underline{k}}.
\label{fromGoursat}
\end{equation}
This bijection is purely algebraic in nature and valid for all $\nu$.

 Lifting invariants give a map
$\cF_{h^k}/\Br_\nu \rightarrow H_2(G^{[k]},C^k,\nu)$.
For any fixed $k$, the Conway-Parker theorem says that this
map is asymptotically a bijection.   Taking the quotient by $\Out(G,C)^{[k]}$ and incorporating
the Goursat conclusion \eqref{fromGoursat} we get 
the desired description of braid orbits:
\begin{equation}
\label{mainbijection1}
\cF_h^{*\underline{k}}/\Br_\nu \stackrel{a \sim}{\longrightarrow} H_2 (G^{[k]},C^k,\nu)/\Out(G,C)^{[k]}.
\end{equation}
The map of \eqref{mainbijection1} is defined
for all allowed $\nu$ and, as indicated by the notation $a \!  \sim$,  is asymptotically a bijection. 

There is, of course, a map $\cF_h^{*\underline{k}}/\Br_\nu \rightarrow \left( \cF_h^*/\Br_\nu \right)^k$;
on the right-hand side of \eqref{mainbijection1}, this corresponds to the natural map 
\begin{equation} \label{mainbijection2}  H_2 (G^{[k]},C^k,\nu)/\Out(G,C)^{[k]} \rightarrow \left(H_2(G, C, \nu)/\Out(G,C)\right)^k.\end{equation} 
Note that the action of $\Out(G,C)^{[k]}$ on $H_2(G^{[k]}, C^k, \nu)$ factors,
under the coordinate projection $H_2(G^{[k]}, C^k, \nu) \rightarrow H_2(G, C, \nu)$, 
through the corresponding    coordinate projection $\Out(G, C)^{[k]} \rightarrow \Out(G,C)$.

\subsection{End of the proof of $\mbox{I} \Rightarrow \mbox{II}$ in the split-cyclic case}  
\label{Endsc} We now assume not only Conditions~1 and 2 of I, but also
Condition $3$.   In this subsection,
we complete the proof of $\mbox{I} \Rightarrow \mbox{II}$ under the auxiliary assumption
that the surjection $G \rightarrow G^{\rm ab}$ is split and $G^{\rm ab}$ is cyclic. 
Some of the notions introduced here are used again in the \S\ref{Endgeneral}, where we complete
the proof without auxiliary assumptions.

Consider  the canonical surjections $H_2(G^{[k]},C^k,\nu) \twoheadrightarrow H_2(G,C,\nu)^k$.   
Under our auxiliary assumption that $G$ has split-cyclic type, Condition 3 and Proposition~\ref{homprop}
show that
$$|H_2(G^{[k]},C^k)| = |H_2(G,C)|^k$$ for all $k$.    Thus, since cardinality does not change when
one passes from groups to torsors, the surjections are bijections.
Moreover, because inner automorphisms act trivially on $H_2(G, C, \nu)$, 
the action of $\Out(G,C)^{[k]}$ on $H_2(G,C,\nu)^k$ actually factors through $\Out(G,C)^k$. 

   Taking the
quotient by $\Out(G,C)^{[k]}$, we can rewrite
  \eqref{mainbijection1} as
\begin{equation}
\label{mainbijection2}
\cF_h^{*\underline{k}}/\Br_\nu \stackrel{a \sim}{\longrightarrow} H_2^*(G,C,\nu)^k.  
\end{equation}
Then standard group theory  shows that the action of $\Br_\nu$ on $\cF^*_h$ is quasi-full
for sufficiently large $\min_i \nu_i$:

In general, consider a permutation group  $B \subseteq \Sym(F)$ with orbit decomposition $F = \coprod_{i=1}^s F_i$.     
Suppose each orbit $F_i$ has size at least $k$.  
  Then the induced action 
of $B$ on $F^{\underline{k}}$ has at least $s^k$ orbits.  If equality holds, then the images $B_i \subseteq \Sym(F_i)$ of
$B$ are each individually $k$-transitive.    If $k \geq 6$, then the classification of finite simple groups
says that $B_i$ contains $\Alt(F_i)$.   Still assuming that $B$ has exactly 
$s^k$ orbits on $F^{\underline{k}}$, it is then elementary that  
$B$ contains $\Alt(F_1) \times \cdots \times \Alt(F_s)$.

\subsection{A lemma on 2-transitive groups}  
\label{lemma2trans}
For the general case, Condition~3 gives 
us control over $\Br_\nu$-orbits only on pairs $(x_1,x_2)$ of distinct elements
in $\cF_h^*$, not tuples of larger length.   To deal with this problem, we replace the
classification of multiply-transitive groups by a statement derived
from the classification of $2$-transitive groups.   The exact formulation
of our lemma is inessential; its import is that full groups are
clearly separated out from other $2$-transitive groups in 
a way sufficient for our purpose.

  \begin{Lemma}  
  \label{translemma} Fix an odd integer $j \geq 5$.    Suppose a $2$-transitive group 
 $\Gamma \subseteq \Sym(X)$ satisfies
 $|X^{\underline{2j}}/\Gamma| \leq 2^{j^2-4j}$.
If  $|X|$ is sufficiently large, then $\Gamma$ is full.
\end{Lemma}

\proof 
To prove the statement, we use the classification of non-full $2$-transitive groups, 
as presented in \cite[\S7.7]{DM}, thereby breaking into a finite number of cases.  For fixed $j$, we discard in each case
a finite number of $\Gamma$ and otherwise establish $|X^{\underline{2j}}/\Gamma| > 2^{j^2-4j}$.

     It suffices to restrict attention to 
 maximal non-full $2$-transitive groups $\Gamma$.  Besides a small number of 
 examples involving seven of the sporadic groups \cite[p.252-253]{DM}, every such maximal $\Gamma$ occurs on the following table.
\[
\begin{array}{rlccc}
\# & \mbox{Type} &\Gamma&  \mbox{Degree $N$} & \mbox{Order $|\Gamma|$} \\
\hline
1 & \mbox{Affine} &AGL_d(p) & p^d & \\
2 & \mbox{Projective} &P\Gamma L_d(q) &  (q^d-1)/(q-1) &  \\ 
3 & \mbox{OS2} & O_{2d+1}(2) &   2^d (2^d \pm 1)/2 & \\ 
4 & \mbox{Unitary} & U_3(q) &   q^3 + 1 & q^3 (q^2-1) (q^3+1)  \\
5 & \mbox{Suzuki} & Sz(q) &  q^2+1 & (q^2+1) q^2 (q-1)  \\
6 & \mbox{Ree}      & R(q)    &  q^3+1 & (q^3+1) q^3 (q-1) \\
\end{array}
\]
The six series are listed in the order they are treated in \cite[p.244-252]{DM}.  
Throughout, $p$ is a prime number and $q=p^e$ is a prime power. 
These numbers are arbitrary, except in Cases 
$5$ and $6$ where the base is $p=2$ and $p=3$ respectively and 
the exponent $e$ is odd.  The orders $|\Gamma|$ in Cases 1-3 
are not needed in our argument and so are omitted from the table.     

   {\em Cases  4-6.} In these cases, the order $|\Gamma|$ grows only polynomially in the degree $N$, with 
$|\Gamma| < N^3$ holding always.    One has
\[
|X^{\underline{2j}}/\Gamma| \geq N^{\underline{2j}}/|\Gamma| > N^{\underline{2j}}/{N^3}.
\]
For $j \geq 5$ fixed and $N \rightarrow \infty$, the right side tends to $\infty$.  So, with
finitely many exceptions, $|X^{\underline{2j}}/\Gamma| > 2^{j^2-4j}$.  

 {\em Case 1.}  In this case,  the affine general linear group $AGL_d(\F_p)$ acts on the affine space 
$\F_p^d$.  Let $w = \min(j,d+1)$.  Fix $x_1$, \dots, $x_w$ in $\F_p^d$ spanning an affine subspace $A$ of 
dimension $w-1$.  The set $A-\{x_1,\dots,x_w\}$ has $p^{w-1}-w$ elements.   There 
are $(p^{w-1}-w)^{\underline{2j-w}}$ ways to 
successively choose $x_{w+1}$, \dots, $x_{2j}$ in $A$ so that all the
$x_i$ are distinct.    The tuples $(x_1,\dots,x_{2j}) \in (\F_p^d)^{\underline{2j}}$ so obtained 
are in different $AGL_d(\F_p)$ orbits.  Thus 
\[
|(\F^d_p)^{\underline{2j}}/AGL_d(\F_p)| \geq \left(p^{w-1}- w \right)^{\underline{2j-w}}.
\]
For fixed $d<j$, so that $w = d+1$, the right side tends to $\infty$ with $p$, and so with finitely many
exceptions, $|(\F^d_p)^{\underline{2j}}/AGL_d(\F_p)|> 2^{j^2-4j}$.  For $d \geq j$, so that $w = j$, one gets no exceptions, as
\[
\left(p^{w-1}- w \right)^{\underline{2j-w}} = \left(p^{j-1}- j \right)^{\underline{j}}  \geq  \left(2^{j-1}- j \right)^{\underline{j}} 
\geq   \left(2^{j-1}- 2j + 1 \right)^{j} > 2^{j^2-4j}.  
\]
[Case 1 is the only case where there is a complicated list of non-maximal $2$-transitive groups.  
Some large ones 
are $AGL_{d/e}(\F_{p^e})  \subset A \Gamma L_{d/e}(\F_{p^e}) \subset AGL_d(p)$, for 
any $e$ properly dividing $d$.]

     Cases $2$ and $3$ are very similar to Case~1, but sufficiently different to require separate
treatments.

{\em Case 2.}  Here $\Gamma = P\Gamma L_d(\F_q) = PGL_d(\F_q).\Gal(\F_q/\F_p)$ 
acts on the projective space $X=\bbP^{d-1}(\F_q)$.  Again let $w = \min(j,d+1)$.  Fix
$x_1$, \dots, $x_w$ in $\bbP^{d-1}(\F_q)$ spanning a projective subspace $P$ of dimension $w-1$.   
Similarly to Case 1, there are $((q^{w}-1)/(q-1) - w)^{\underline{2j-w}}$
ways to successively choose $x_{w+1}$, \dots, $x_{2j}$ in $P$ so that all the
$x_i$ are distinct.    The tuples $(x_1,\dots,x_{2j}) \in \bbP^{d-1}(\F_q)^{\underline{2j}}$ so obtained 
are in different $PGL_d(\F_q)$ orbits.   However one $P\Gamma L_d(\F_q)$ orbit can consist of 
up to $e$ different $PGL_d(\F_q)$ orbits.   Thus our lower bound in this case is
\[
|\bbP^{d-1}(\F_q)^{\underline{2j}}/P\Gamma L_d(\F_q)| \geq \frac{1}{e} \left( \frac{q^{w}-1}{q-1} - w \right)^{\underline{2j-w}}.
\]
Again the subcase $d < j$, where $w = d+1$, is simple:  the right side tends to $\infty$ with $q$ and 
so $|\bbP^{d-1}(\F_q)^{\underline{2j}}/P\Gamma L_d(\F_q)| > 2^{j^2-4j}$ holds
with only finitely many exceptions.  
For $d \geq j$,  so that $w=j$ again, one has no further exceptions as
\[
\frac{1}{e} \left(  \frac{q^{w}-1}{q-1}- w \right)^{\underline{2j-w}}
> \frac{1}{e} \left(q^{j-1} - 2j+1 \right)^{j} >   \left(2^{j-1} - 2j+1 \right)^{j} > 2^{j^2-4j}.
\]

{\em Case 3.}  Here the group in question in its most familiar guise is $\Gamma = Sp_{2d}(\F_2)$
for $d \geq 2$.
It is better in our context to view $\Gamma = O_{2d+1}(\F_2)$, as from
this point of view the 2-transitive actions appear most naturally.  In fact 
the orbit decomposition of the natural action of $O_{2d+1}(\F_2)$ is 
\[
\F_2^{2d+1} - \{0\} = X_{-1} \coprod X_1 \coprod X_0.
\]
Here $X_0$ is the set of isotropic vectors.  The pair 
$(O_{2d+1}(\F_{2}),X_0)$ is a copy of the more standard
pair $(Sp_{2d}(\F_2),\F_2^{2d}-\{0\})$ and so in
particular $|X_0|=2^{2d}-1$.  A non-isotropic vector
is in $X_1$ if its stabilizer is the split orthogonal group
$O^+_{2d}(\F_2)$ and is in $X_{-1}$ if its stabilizer is the non-split
orthogonal group $O^-_{2d}(\F_2)$.  From the order of the stabilizers
one gets that $|X_\epsilon| = 2^{d-1}(2^d+\epsilon)$.  
  While the action 
of $\Gamma$ on $X_0^{\underline{2}}$ has two orbits,
the actions on the other two $X_\epsilon$ are $2$-transitive.  
[Familiar examples for $O_{2d+1}(\F_2) = Sp_{2d}(\F_2)$ come from $d=2$, and $d=3$.  
Here the groups respectively are $S_6$, and $W(E_7)$.   The
orbit sizes on $(X_{-1},X_{1},X_0)$ are  $(6,10,15)$ and $(28,36,63)$ respectively.]

By discarding a finite number of $\Gamma$, we
 can assume $d \geq j$.  Fix $x_1$, \dots,  $x_j$ in $X_\epsilon$ spanning a $j$-dimensional
vector space $V \subset \F_2^{2d+1}$ on which the quadratic form remains non-degenerate and each
$x_i$ has type $\epsilon$ in this smaller space.   Let 
$V_\epsilon = V \cap X_\epsilon$.  Writing $j = 2u+1$, one has $|V_\epsilon| = 2^{u-1} (2^u + \epsilon)$.
There are $(|V_\epsilon|-j)^{\underline{j}}$ ways to successively choose 
$x_{j+1}$, \dots, $x_{2j}$ in $V_\epsilon$ so that all the
$x_i$ are distinct.  One has
\[
|X_\epsilon^{\underline{2j}}/O_{2d+1}(\F_2)| \geq \left(2^{u-1} (2^u+\epsilon) -j \right)^{\underline{j}}
  \geq \left(2^{u-1} (2^u+\epsilon) -2j+1 \right)^{{j}} > 2^{j^2-4j}.
\]
Thus there are no further exceptional $\Gamma$ from this case. \qed

\subsection{End of the proof if I $\Rightarrow$ II  in general}  
\label{Endgeneral}
We now end the proof without the split-cyclicity assumption,
by modifying the standard argument of \S\ref{Endsc}.

Consider again the diagram \eqref{ei2} relating two five-term exact sequences.
The last three terms of the top sequence
and the last four terms of the bottom sequence give respectively 
\begin{eqnarray*}
|H_2(G^{[k]} )| & \leq & |H_2(G' )_{G^{\rm ab}}|^k   |H_2(G^{\rm ab})|, \\
|H_2(G' )_{G^{\rm ab}}|^k & \leq &  \frac{ |H_3(G^{\rm ab} )|^k |H_2(G )|^k}{|H_2(G^{\rm ab} )|^k}.
\end{eqnarray*}
 Combining these inequalities and replacing
$H_2(G^{[k]} )$ by its quotient $H_2(G^{[k]},C^k)$ yields
\begin{equation} \label{bound}  | H_2(G^{[k]}, C^k) | \leq  | H_2(G  ) \times H_3(G^{\ab} )|^k. 
\end{equation}
As described in \S\ref{Endsc}, Condition~3 implies that for $\min \nu_i$ 
sufficiently large, the action of $\Br_\nu$ on $\cF_h^*$ is $2$-transitive
when restricted to each orbit. 
We will use this $2$-transitivity and the exponential bound \eqref{bound} to conclude that
the action of $\Br_\nu$ on $\cF^*_h$ is asymptotically quasi-full.  

Consider $S_m$ in its standard full action on $Y_m = \{1,\dots,m\}$.   
The induced action on $X_m = Y_m \coprod Y_m$ is not quasi-full.   Let
$a_{k,m}$ be the number of orbits of $S_m$ on $Y_m^{\underline{k}}$.
As $m$ increases the sequence $a_{k,m}$ stabilizes at a number
$a_k$.  The sequence $a_k$ appears in \cite{OEIS} as A000898.  
There are several explicit formulas and combinatorial interpretations.
The only important thing for us is that $a_k$ grows superexponentially,
as indeed $a_k/a_{k-1} \sim \sqrt{2 k}$.    

From \eqref{bound} we know that there exists an odd number $j$ with 
\[
|H_2(G^{[2j]}, C^{2j},\nu)/\Out(G, C)^{[2j]}| \leq | H_2(G^{[2j]}, C^{2j}) | < \min(2^{j^2-4j},a_{2j}).
\]
 By \eqref{mainbijection1}, the left-hand set  is identified with 
$|\cF_h^{* \underline 2j}/\Br_\nu|$ for sufficiently large 
$\min_i \nu_i$.   Lemma~\ref{translemma} above says that, at the possible 
expense of
making $\min_i \nu_i$ even larger, each  orbit
of the action of $\Br_\nu$ on $\cF_h^*$ is full.   
Our discussion of the action of $S_m$ on $Y_m$ 
says that the constituents are pairwise 
non-isomorphic, again for sufficiently large $\min_i \nu_i$.   The classical Goursat lemma then says 
the action is quasi-full.   \qed  

A consequence of the results of the section is that in fact the 
equivalence 3 $\Leftrightarrow$ $\hat{3}$ of Proposition~\ref{homprop} holds without
the assumption of split-cyclicity.   Condition~E is also meaningful
in general, and it would be interesting to identify the class
of $(G,C)$ for which the equivalence extends to include E.

 \section{Proof of II $\Rightarrow$ I}
\label{reverse}
    In this section, we complete the proof of Theorem~\ref{mt} by proving that (not I) implies (not II). 
    Accordingly, we fix a centerless group $G$ and a list $C = (C_1,\dots,C_r)$ of conjugacy classes
     and consider consequences of the failure
     of Conditions 1, 2, and 3 in turn.   In all three cases, we show more than is needed for 
     Theorem~\ref{mt}.

 \subsection{Failure of Condition 1} 
 \label{notpseudosimple}  The failure of the first condition 
 requires a somewhat lengthy analysis, because it breaks into two quite different 
 cases.   The conclusion of the following lemma shows more than
 asymptotic quasi-fullness of $\AHur_h^* \rightarrow \AConf_\nu$ fails; it shows
 that asymptotically each individual component $\AHur_{h,\ell}^* \rightarrow \AConf_\nu$ 
 fails to be full.  
 
 \begin{Lemma}
 \label{fail1lem} Let $G$ be a centerless group which is not pseudosimple. Let $C = (C_1,\dots,C_r)$ be 
 a list of conjugacy classes.  Consider varying allowed $\nu \in \Z_{\geq 1}^r$ and
 thus varying Hurwitz parameters $h = (G,C,\nu)$.  Then for 
 $\min_i \nu_i$ sufficiently large and any $\ell \in H_h^*$, the action of 
 $\Br_\nu$ on 
 $\cF_{h,\ell}^*$ is not full.    
 \end{Lemma}
\proof A group is pseudosimple exactly when it satisfies two conditions:  
 $A$, it has no proper nonabelian quotients; 
$B$, its derived group is nonabelian.  
We assume first that $A$ fails.  Then we assume that $A$ holds but $B$ fails. 
\medskip

\noindent{\em Assume $A$ fails.}
   Let $\overline{G}$ be a proper nonabelian quotient.  Let $\bar{h} = (\overline{G},(\overline{C}_1,\dots,\overline{C}_r),\nu)$ 
 be the corresponding quotient Hurwitz parameter.   Consider the natural map 
 $H_h \rightarrow H_{\bar{h}}$  from \S \ref{funcmap} and
 let $\bar{\ell}$ be the image of $\ell$.
 
  By the definition of Hurwitz parameter, the classes $C_i$ generate $G$.  At least one of the surjections 
 $C_i \rightarrow \overline{C}_i$ has to be non-injective, as otherwise the kernel of $G \rightarrow \overline{G}$
 would be central in $G$ and  $G$ is centerless.  
 So $|C_i| \geq 2 |\overline{C}_i|$ for at least one $i$.    Similarly, since $\overline{G}$ is nonabelian and 
 generated by the $\overline{C}_i$, one has $|\overline{C}_i| \geq 2$ for at least one $i$.   
 
We now examine
the induced map 
$\cG_{h, \ell} \rightarrow \cG_{\bar{h}, \bar{\ell}}$.  Let 
 $\cI_{h,\ell}$ be its image and $\phi_{h,\ell}$ be the size of its largest fiber.
 We will use the two inequalities of the previous paragraph to show that both 
$\phi_{h,\ell}$ and $|\cI_{h,\ell}|$ grow without bound with $\min_i \nu_i$.

From $|C_i| \geq 2 |\overline{C}_i|$ and two applications of the asymptotic mass formula \eqref{ff8},
 one gets $|\cG_{h,\ell}| \geq 1.5^{\min_i \! \nu_i} |\cG_{\bar{h},\bar{\ell}}|$
and hence $\phi_{h,\ell} \geq 1.5^{\min_i \! \nu_i}$. 

To show the growth of $|\cI_{h,\ell}|$, we assume without loss of generality that $|\overline{C}_1| \geq 2$ and
choose $y_1 \neq y_2 \in \overline{C}_1$. 
Let $M$ be the exponent of a reduced Schur cover $\tilde{G}_C$ of $G$.  Let $k$ be a positive integer  and let 
$a_1$, \dots, $a_k$ be a sequence with $a_i \in \{1,2\}$.  Then for  $\min_i \nu_i$ large enough,
we claim that $\cI_{h,\ell}$ contains an element of the form
\begin{equation} \label{Fred}(\underbrace{y_{a_1}, \dots, y_{a_1}}_M,  \dots, \underbrace{y_{a_k}, \dots, y_{a_k}}_M,
\underbrace{x_1, \dots, x_{\nu_1-Mk}}_{{\rm all} \; {\rm in} \; \overline{C}_1}, \dots,
\underbrace{x_{n-kM-\nu_r+1}, \dots, x_{n-kM}}_{{\rm all} \; {\rm in} \; \overline{C}_r} ). \!\!\! \end{equation} 
To see the existence of such an element, fix a lift  $C_i^*$ of the conjugacy class $C_i$ to $\tilde{G}_C$
and choose $\tilde{y}_1, \tilde{y}_2 \in C_1^*$ mapping (under $\tilde{G}_C \rightarrow G \rightarrow \bar{G}$) to $y_1,y_2 \in \overline{C}_1$ respectively. 
 
Let $z \in H_2(G,C)$ be chosen so that $z^{-1} \cdot \prod_{i} [C_i]^{\nu_i} = \ell$ inside $H_2(G, C, \nu)$. Consider the equation 
\begin{equation}
\label{Fred2}
(\tilde{y}_{a_1}^M  \cdots \tilde{y}_{a_k}^M) \underbrace{ \tilde{x}_1 \cdots  \tilde{x}_{\nu_1-kM}}_{{\rm all} \; {\rm in} \;C_1^*} \cdots  \underbrace{\tilde{x}_{n-kM-\nu_r+1}\cdots \tilde{x}_{n-kM}}_{{\rm all} \; {\rm in} \; C_r^*}  = z,
\end{equation}
where $\tilde{x}_i \in C_i^*$.     By our choice of $M$, the powers $\tilde{y}^M_{a_i}$ are all 
the identity in $\tilde{G}_C$.   One has $[C_1^*]^{\nu_1-kM} \cdots [C_r^*]^{\nu_r} = [z]$ in $\tilde{G}_C^{\rm ab} = G^{\rm ab}$,
both sides being the identity.  The asymptotic mass formula then applies to say that \eqref{Fred2} in fact has many solutions 
$(\tilde{x}_1,\dots,\tilde{x}_{n-kM})$ where moreover $\tilde{x}_i$ generate $\tilde{G}_C$. 
Now, the image of $(\tilde{y}_{a_1}, \dots, \tilde{y}_{a_k}, \tilde{x}_1, \dots, \tilde{x}_{n-kM})$ actually
defines an element of $\cG_{h,\ell}$, and its 
image in $\bar{G}$ is an element of $\cI_{{h},{\ell}}$ of the
form \eqref{Fred}.   Varying $(a_1,\dots,a_k)$ now, always taking $\min_i \nu_i$ sufficiently large, we conclude $|\cI_{h,\ell}| \geq 2^k$.  

For large enough $\min_i \nu_i$ the action of $\Br_{\nu}$ on $\cG_{h, \ell}$ is transitive, by the Conway--Parker theorem. 
This action preserves a partition of $\cG_{h,\ell}$  
into $b = |\cI_{h,\ell}|$ blocks, each of size $f = \phi_{h,\ell}$.    Thus the image of $\Br_\nu  $ on $\cG_{h,\ell}$ 
  is contained in the wreath product $S_f \wr S_b$.   Hence the image of 
  $\Br_{\nu}$ on $\cF_{h,\ell}^*$ is contained in a subquotient of $S_f \wr S_b$.
  But we have established that $f$ and $b$ increase indefinitely
  with $\min_i \nu_i$.   Let $a = |\Aut(G,C)_\ell|$ and $m=|\cF^*_{h,\ell}|$ so that 
  $|\cG_{h,\ell}| = ma = fb$.  
  As soon as $\min(f,b)>a$, one has $m > \max(f,b)$ and the alternating group $A_m$ is not a subquotient
  of $S_f \wr S_b$.  So the action of $\Br_\nu$ on $\cF^*_{h,\ell}$ is not full.

\medskip
   
\noindent {\em Assume $A$ holds but $B$ fails.}  The assumptions force $G'$ to be isomorphic
to the  additive group of $\F_p^w$ for some prime $p$ and some power $w$.  Moreover, consider the action
 of $G^{\rm ab}$ on $G'$.   Now  $G'$, considered as an $\F_p$-vector space, is an  irreducible representation of $\F_p[G^{\rm ab}]$.
 $G^{\rm ab}$ must have order coprime to $p$, as otherwise the fixed subspace for the $p$-primary part of 
 $G^{\rm ab}$   would be
 a proper subrepresentation.     So $\F_p[G^{\rm ab}]$
 is isomorphic to a sum of finite fields and the action on $G' = \F_p^w$ is through a 
 single summand $\F_q$.   We can thus identify $G'$ with  the additive group of a finite field $\F_q$ and $G^{\rm ab}$ with a subgroup
 of $\F_q^{\times}$ in such a way that $G$ itself is a subgroup of the affine group 
 $\F_q.\F_{q}^\times$.    Moreover, $G^{\rm ab} \subseteq \F_q^{\times}$ acts irreducibly on $\F_q$
 as an $\F_p$-vector space. 
 
 We think of elements of $G$ as affine transformations $x \mapsto mx+b$.  
 Since braid groups act on the right in \eqref{bdef}, we compose these affine transformation from
 left to right, so that the group law is
 $\binom{m_1}{b_1} \binom{m_2}{b_2} = \binom{m_1m_2}{m_2b_1+b_2}$.  

We think of elements $(g_1,\dots,g_n) \in \cG_h$ with $g_i = \binom{m_i}{b_i}$ in terms of the following matrix:
\begin{equation}
\label{gmat}
\left[
\begin{array}{cccccc}
m_1& \cdots & m_i & m_{i+1}& \cdots & m_n  \\
b_1 & \cdots & b_i  &  b_{i+1} &  \cdots & b_n 
\end{array}
\right]
\end{equation}
The top row is determined by $C$, via $m_i = [C_i]$.  
Thus, via the bottom row, we have realized $\cG_h$ as a subset of $\F_q^n$.    
We can assume without loss of generality that none of the $C_i$ are the identity class.  
Then the requirement $g_i \in C_i$ for membership in $\cG_h$ gives $|G^{\rm ab}|$ choices 
for $b_i$ if $m_i=1$.  If $m_i \neq 1$ then $g_i \in C_i$ allows all $q$ choices for $b_i$.

Now briefly view $(g_1,\dots,g_n)$ as part of the larger catch-all set $G^n$ of \S\ref{catchall}, on which
the standard braid operators $\sigma_i$ act.    The braiding rule  \eqref{bdef} in our current setting becomes
{\small
\[
\left( \dots, \binom{m_i}{b_i}, \binom{m_{i+1}}{b_{i+1}}, \dots \right)^{\sigma_i} = 
\left(\dots,\binom{m_{i+1}}{b_{i+1}}, \binom{m_i}{b_{i+1} + m_{i+1} b_i - m_i b_{i+1}}, \dots \right).
\]
}

\noindent 
Thus the action of $\sigma_i$ corresponds to the the bottom row of \eqref{gmat},
viewed as row vector of length $n$, being multiplied on the right by an $n$-by-$n$ matrix 
in $GL_n(\F_q)$.    

Returning now to the set $\cG_h$ itself, any element of $\Br_\nu$ can be written as a product of 
the $\sigma_i$ and their inverses.  Accordingly, image of the $\Br_\nu$ in  $\Sym(\cG_h)$
lies in  $GL_n(\F_q)$.

To prove non-fullness, it
suffices to bound the sizes of groups.
On the one hand, 
\[
|\mbox{Image of $\Br_\nu$ in  $\Sym(\cF_{h,\ell}^*)$}| \leq |\mbox{Image of $\Br_\nu$ in  $\Sym(\cG_h)$}| \leq  |GL_n(\F_q)| < q^{n^2}.
\] 
 On the other hand, let $b = |H_2(G,C)| |\Out(G,C)|+1$.  Then,
 using  \eqref{ff8},
\eqref{refined-mass} and the fact that $|C_i| \in \{|G^{\rm ab}|,q\}$, 
one has 
\[
|\cF^*_{h,\ell}| >  \frac{\prod_i |C_i|^{\nu_i}}{|G| |G'| b} \geq \frac{|G^{\rm ab}|^{n-3}}{q^2 b},
\]
for all sufficiently large $n$.  
Certainly
$
q^{n^2} < \frac{1}{2} \left(\frac{a^{n-3}}{q^2 b}\right)!
$
for any fixed $a$, $b$, $q>1$ and sufficiently large $n$.   Thus  
the image of $\Br_\nu$ in $\Sym(\cF_{h,\ell}^*)$ cannot contain $\Alt(\cF^*_{h,\ell})$.  
  \qed
  
  The paper \cite{EEHS} calculates monodromy in cases with $G = S_3$ and $G=S_4$, providing worked out
examples.  Another illustration of the case with affine monodromy is  \cite[Prop.\ 10.4]{MM}.

\subsection{Failure of Condition 2.}  \label{ambiguous}  Our next lemma has the  same conclusion as the previous lemma.  

 \begin{Lemma}
 \label{fail2lem} Let $G$ be a centerless group.  Let $C = (C_1,\dots,C_r)$ be 
 a list of conjugacy classes with at least one $C_i$ ambiguous.  Consider varying allowed $\nu \in \Z_{\geq 1}^r$ and
 thus varying Hurwitz parameters $h = (G,C,\nu)$.  Then for 
 $\min_i \nu_i$ sufficiently large and any $\ell \in H_h^*$, the action of 
 $\Br_\nu$ on 
 $\cF_{h,\ell}^*$ is not full.    
 \end{Lemma}

\proof Introduce indexing sets $B_i$ by 
writing
\[
C_i = \coprod_{b \in B_i} C_{ib},
\]
where each $C_{ib}$ is a single $G'$ orbit.   Our hypothesis says that 
at least one of the $B_i$  -- without loss of generality, $B_1$ -- has size larger than $1$. On the other hand,
at least one of the $B_i$ has size strictly less than $C_i$; otherwise
$G'$ would centralize each element of each $C_i$, and then all of $G$, which
is impossible for $G$ center-free.    

Define
\[
\cG_h^{\rm amb} = \overbrace{B_1 \times \cdots \times  B_1}^{\nu_1} \times \cdots \times  \overbrace{B_r \times \cdots \times  B_r}^{\nu_r}.
\]
The group $G$ acts transitively through its abelianization $G^{\rm ab}$ on each $B_i$. 
For a lifting invariant $\ell \in H_h$, consider the natural map
$\cG_{h,\ell} \rightarrow \cG_h^{\rm amb}$.   The action of the braid group $\Br_{\nu}$ on $\cG_{h,\ell}$
descends to an action on $\cG_h^{\rm amb}$.  

Now we let $\min_i \nu_i \rightarrow \infty$ and get the following consequences, by 
arguments very closely paralleling those for the first case of Lemma~\ref{fail1lem}.
First, the image of the map $\cG_{h,\ell} \rightarrow \cG_h^{\rm amb}$,
 has size that goes to $\infty$.    Second, the mass formula again shows
that $\frac{|\cG_{h,\ell}|}{|\cG_h^{\rm amb}|} \rightarrow \infty$ with $\min_i \nu_i$.  By the last paragraph
of the first case of the proof of Lemma~\ref{fail1lem},  the action of $\Br_\nu$ on each orbit of $\cF_{h,\ell}^*$ is
 forced to be imprimitive, and hence not full. \qed

For a contrasting pair of examples, consider 
 $h = (S_5,(C_{2111},C_{311},C_5),\nu)$ for 
$\nu = (2,2,1)$ and $\nu=(2,1,2)$.   The monodromy group for the former is all of $S_{125}$, 
despite the presence of the ambiguous class $C_5$.    The monodromy group for the latter
is $S_{85} \wr S_2$ and represents the asymptotically-forced non-fullness.

\subsection{Failure of Condition 3}  The last lemma of this section is different in structure 
from the previous two, and its proof is essentially a collection of some of our previous arguments.   From the discussion of surjectivity after \eqref{ei2}, one always has
\begin{equation}
\label{defa}
|H_2(G^{[2]},C^2)| = a |H_2(G,C)|^2
\end{equation}
for some positive integer $a$.   Condition~3 is that $a=1$. The number 
$a$ reappears as the cardinality of every fiber of 
the maps of torsors 
$$H_2(G^{[2]}, C^2, \nu)  \stackrel{\pi}{\twoheadrightarrow} H_2(G,C,\nu)^2$$
considered in \S\ref{funcmap}.

Now suppose that $\nu$ is such that all $\nu_i$ are divisible
by both the exponent of $H_2(G,C)$ and the exponent of $H_2(G^{[2]}, C^2)$.    In that case, 
we have identifications 
\begin{equation}
  \begin{array}{ccc}
   H_2(G^{[2]}, C^2, \nu) &  \stackrel{\pi}{\rightarrow}& H_2(G,C,\nu)^2   \\
  f \rotatebox{270}{$\!\!\!\! \rightarrow$} & & g\rotatebox{270}{$\!\!\!\! \rightarrow$}    \\
    H_2(G^{[2]}, C^2) & \rightarrow & H_2(G,C)^2 
  \end{array}
  \end{equation} 
where the vertical  bijections $f,g$ come from  \S \ref{subsec:torsors}, and the
fact that the diagram commutes is also explained there. 
 
Then $E= \pi^{-1} g^{-1}(0) \subseteq  H_2(G^{[2]}, C^2, \nu)$
is a fiber of $\pi$. It has size $\geq 2$ and $f(E)   \subseteq H_2(G^{[2]}, C^2)$ is a subgroup. 

The group $\Out(G, C)^{[2]}$, defined before \eqref{fromGoursat}, acts on $H_2(G^{[2]}, C^2, \nu)$
and also (compatibly) on $H_2(G^{[2]}, C^2)$. 
It preserves $E$ and acts on it with at least two orbits, because
it fixes the zero element of $f(E)$.  Under the bijection \eqref{mainbijection1}, 
these two orbits correspond to two different  braid orbits $O, O'$ on
 $\left( \mathcal{F}_h^* \right)^{\underline{2}}$
  which project (in both coordinates) to the same braid orbit on $\mathcal{F}_h^*$.

To summarize, we have proved:

\label{3fail}   
\begin{Lemma}
 Let $G$ be a pseudosimple group, let $C = (C_1,\dots,C_r)$ be 
 a list of unambiguous conjugacy classes, and 
 suppose  $a>1$ in \eqref{defa}.

  Consider $\nu$ with 
 all $\nu_i$ a multiple of the exponent of both $H_2(G,C)$ and $H_2(G^{[2]}, C^2)$.   
  Identify $H_2(G,C,\nu) = H_2(G,C)$ as in \S\ref{subsec:torsors}, writing $0 \in 
 H_2(G,C,\nu)$ for the  element corresponding to $0$. Similarly identify $H_2(G^{[2]},C^2,\nu) = H_2(G^{[2]},C^2)$.  

Then for
$\min_i \nu_i$ sufficiently large, the action of $\Br_\nu$ on 
 $\cF_{h,0}^*$ is not $2$-transitive and hence not full.  
 \end{Lemma}

 \section{Full number fields} \label{FullProposal} 
 
        Theorem~\ref{mt} guarantees the existence of infinitely many  quasi-full covers 
           $\pi_{h}^* : \AHur_{h}^* \rightarrow  \AConf_\nu$
       associated to each simple group $T$.     As discussed in \S\ref{descent}, if all the $C_i$ are different and 
       conjugate classes $C_i$ occur with 
       equal multiplicity, then $\pi_h^*$ canonically descends to a covering of $\Q$-varieties,
       \begin{equation} \label{pihl} \pi_{h}^* : \Hur_{h}^* \rightarrow  \Conf^\rho_\nu. \end{equation}
       This final section explains why we expect specializations of these
  covers to give enough fields for Conjecture~\ref{mc}.   
  
   Our object here is to give an
  overview only, as we defer a more detailed treatment to \cite{HNF}.  In particular, we return
  to the setting of \S\ref{introNF}, considering only $h$ where all $C_i$ are individually rational.
  Then the twisting $\rho$ is trivial and the base of \eqref{pihl}  is just $\Conf_{\nu}$,
  as defined in \S\ref{configspace}.

 \subsection{Specialization}   First, we give a few more details on the
 specialization process. 
 The $\Q$-variety $\Conf_{\nu}$ has a natural structure of scheme over $\Z$.  In particular,
one says that a point $u \in \Conf_{\nu}(\Q)$ is $\primes$-integral
if  it belongs to $\Conf_{\nu}(\Zprimes)$.   
Concretely, a $\primes$-integral point $u$ can be specified by giving 
  binary homogeneous forms $(q_1, \dots, q_r)$, where
   \begin{equation} \label{binary}   \mbox{$q_i \in \Z[x,y]$ and  $\disc(\prod q_i)$ is divisible only by primes in $\primes$.}
   \end{equation}
  To avoid obtaining duplicate fields in the specialization process, one can normalize in 
  various ways to take one point from each $PGL_2(\Q)$ orbit intersecting 
  $\Conf_\nu(\Zprimes)$.  This is done systematically in \cite{RPoly} and 
  these sets of representatives are arbitrarily large for any given non-empty
   $\primes$.   
 
 \subsection{Rationality of components}   For $h$ to be useful in supporting Conjecture~\ref{mc}, 
 it is essential that the subcover
 $\AHur_{h,\ell}^* \rightarrow \AConf_\nu$ is defined over $\Q$ for at least one lifting invariant
 $\ell \in H_h^*$.   In this subsection, we 
 explain that for fixed $(G,C)$, many $h = (G,C,\nu)$ may not have
 such a rational $\ell$, but infinitely many do.   
 
Consider the lifting invariant map, 
 \begin{eqnarray*}
\inv^*_h : \pi_0( \AHur^*_h) \rightarrow H^*_h.
 \end{eqnarray*}
 Since $\AHur^*_h = \Hur^*_h(\C)$ is the set of complex points of a $\Q$-variety, there
 is a natural action of $\Gal(\overline{\Q}/\Q)$ on $\pi_0( \AHur^*_h)$.   
 Likewise, via its standard action on conjugacy classes of groups, 
 $\Gal(\overline{\Q}/\Q)$ acts on $H^*_h$.   This
 latter action is through the abelianization $\Gal(\overline{\Q}/\Q)^{\rm ab}$ 
 and can be calculated via character tables.   With these two actions,
 the lifting invariant map is equivariant up to sign (i.e., the action of $\sigma \in \Gal$ on one side
 corresponds to $\sigma^{\pm 1}$ on the other; we didn't compute the sign) -- see  \cite[v1, \S 8]{EVW2}.  
 
 To make the issue at hand more explicit,  
  suppose $|H^*_h|=2$.  
 Then $\Gal(\overline{\Q}/\Q)^{\rm ab}$ acts on $H_h^*$ the same way it acts on
 the complex numbers $\{-\sqrt{d},\sqrt{d}\}$  for some square-free integer $d$.
 The case $d \neq 1$ is common, and then any specialized field $K^*_{h,u}$ contains 
 $\Q(\sqrt{d})$ and hence---outside of the trivial case $K^*_{h,u}=\Q(\sqrt{d})$---is not 
 full.   The case $d=1$ is more favorable to us, as commonly $K^*_{h,u}$ factors 
 into two full fields.    Explicit examples of both $d \neq 1$ and $d=1$ can be easily built 
 \cite{HNF}.
 
 To see in general that for infinitely many $\nu$, the action of $\Gal(\overline{\Q}/\Q)^{\rm ab}$ on
  $H_h^*=H_{G,C,\nu}^*$ has at least one fixed point, we apply 
  the simple remark from  \S\ref{subsec:torsors}.
  Namely suppose that all $\nu_i$ are multiples of the exponent of $H_2(G,C)$.  Then,
the torsor $H_h$ can be canonically
identified with $H_2(G,C)$ itself.   Then $\Gal(\overline{\Q}/\Q)^{\rm ab}$ fixes the identity element
of $H_h$, and so also fixes the image of the identity  in $H_h^*$.

\subsection{A sample cover}  
\label{cover25} To illustrate the ease of producing full fields, we summarize here
the introductory example of \cite{HNF}.    For this example, we take $h = (S_5,(C_{2111},C_5),(4,1))$.
Then $\AHur_h^* = \AHur_h$ is a full cover of $\AConf_{4,1}$ of degree $25$.  
The fiber of 
$\AHur_h \rightarrow \AConf_{4,1}$ over
the configuration $u = (D_1,D_2) = (\{a_1, a_2, a_3, a_4\},\{\infty\})$
consists of all equivalence classes of quintic polynomials
 \begin{equation}
\label{fdef}
g(y) = y^5 + b y^3 + c y^2 + d y + e
\end{equation}
whose critical values are $a_1, a_2, a_3, a_4$.   Here the equivalence
class of $g(y)$ consists of the five polynomials $g(\zeta y)$ where $\zeta$ runs over fifth
roots of unity.  

Explicitly, consider the resultant $r(t)$ of $g(y)-t$ and $g'(y)$.   Then $r(t)$ equals 
{\small \begin{eqnarray*}
 \nonumber \!\!\!\!\!&&3125 t^4 +1250 (3 b c-10 e) t^3 \\ 
\nonumber \!\!\!\!\!&&+\left(  108 b^5-900 b^3 d+825 b^2 c^2-11250 b c e+2000 b d^2+2250 c^2 d+18750    e^2   \right) t^2  \\
\nonumber  \!\!\!\!\!&&  -2 \left( 108 b^5 e-36 b^4 c d+8 b^3 c^3-900 b^3 d e+825 b^2 c^2  e+280 b^2 c d^2 \right.  \\ 
\label{resform} \!\!\!\!\!&& \left. 
    \qquad -315 b c^3 d-5625 b c e^2+2000 b d^2 e+54 c^5+2250 c^2 d e-800 c d^3+6250 e^3
 \right)t \\
\nonumber  \!\!\!\!\!&&  + \left(108 b^5 e^2-72 b^4 c d e+16 b^4 d^3+16 b^3 c^3 e-4 b^3 c^2 d^2-900 b^3 d
    e^2+825 b^2 c^2 e^2  \right. \\ 
    \nonumber \!\!\!\!\!&& \left. \qquad
    +560 b^2 c d^2 e-128 b^2 d^4-630 b c^3 d e+144 b c^2
    d^3-3750 b c e^3 \right. \\ 
    \nonumber 
    \!\!\!\!\!&& \left. \qquad +2000 b d^2 e^2+108 c^5 e-27 c^4 d^2+2250 c^2 d e^2-1600 c
    d^3 e+256 d^5+3125 e^4\right).
\end{eqnarray*}}
\vspace{-.2in}

\noindent  For fixed $\{a_1, a_2, a_3, a_4\}$, there are generically $125$ different solutions $(b,c,d,e)$ to the equation  
$r(t) =  3125 (t-a_1) \cdots (t-a_4)$.  Two solutions are equivalent exactly if they have the same $e$.  Whenever
$D_1$ is rational, i.e.\   $\prod (t-a_i) \in \Q[t]$,  the set of $e$ arising forms the set of
roots of a degree $25$ polynomial with rational coefficients.  By taking 
$u \in \Conf_{4,1}(\Z[\frac{1}{30}])$ one gets more than $10000$ different fields with 
Galois group $A_{25}$ or $S_{25}$ and discriminant of the form
$\pm 2^a 3^b 5^c$.

\subsection{Support for Conjecture~\ref{mc}}

\label{support} 
 Let $F_\cP(m)$ be the number of full fields 
ramified within $\cP$ of degree $m$.  The mass heuristic \cite{B} gives an expected value
 $\mu_\cP(m)$ for $F_\cP(m)$ as an easily computed product of local masses.  
 This heuristic has had clear success in the setting of 
 fixed degree and large discriminant, 
 being for example exactly right on average for $m=5$  \cite{B5}.

The numerical support
 for Conjecture~\ref{mc} presented in \cite{HNF} gives evidence
 that specialization of the covers \eqref{pihl} does indeed
 behave generically.  General computations for fixed $\cP$ in 
 arbitrarily large degree do not seem possible.  However   our numerical support at least shows that 
 specialization of Hurwitz covers produces many fields in degrees
 larger than would be expected from the mass heuristic.

 For instance, one of many examples in \cite{HNF} 
 comes from the Hurwitz parameter  
 $h=(S_6,(C_{321},C_{2111},C_{3111},C_{411}),(2,1,1,1))$. 
 The covering $\Hur_h \rightarrow \Conf_{2,1,1,1}$
is full of  degree $202$.  
 The specialization set 
 $\Conf_{2,1,1,1}(\Z[\frac{1}{30}])$ intersects
 exactly $2947$ different $\PGL_2(\Q)$ orbits 
 on the set $\Conf_{2,1,1,1}(\Q)$ \cite{RPoly}.  
  The mass heuristic predicts $\sum_{m=202}^\infty \mu_{\{2,3,5\}}(m) < 10^{-16}$ full
 fields in degree $\geq 202$. 
 However
 specialization is as generic as it could be, as the $2947$ algebras $K_{h,u}$ are 
pairwise non-isomorphic and  all full. 

Even sharper contradictions to the mass heuristic are obtained in \cite{RCheb}
from fields ramified at just two primes.  However the construction
there is very special, and does not give fields in arbitrarily large degree
for a given $\cP$.    Here we have not just the large supply of full covers
studied in this paper, but also very large specialization sets \cite{RPoly}
giving many opportunities for full fields.    
Specialization in large
degrees would have to behave extremely non-generically 
for Conjecture~\ref{mc} to be false.   
Our belief is that Conjecture~\ref{mc} still holds
with the conclusion strengthened to 
$F_\cP(m)$ being unbounded.

\subsection{Concluding discussion}  
There are other aspects of the sequences $F_{\cP}(m)$ that 
are not addressed by our Conjecture~\ref{mc}.   Most notably, the 
fields arising from full fibers of Hurwitz 
covers occur only in degrees for which there is a cover.  By the 
mass formula, these degrees form a sequence
of density zero.

A fundamental question is thus the support of the sequences
  $F_{\cP}(m)$, meaning
the set of degrees $m$ for which $F_{\cP}(m)$ is positive.  
One extreme 
possibility,  giving as much credence to the mass heuristic
as is still reasonable, is that $F_{\cP}(m)$ has support 
on a sequence of density zero in general and is eventually 
zero unless $\cP$ contains the set of prime divisors of the order of a finite simple group. 
This would imply that the classification of finite simple groups has an unexpected 
governing influence on a part of 
algebraic number theory seemingly
quite removed from general group theory.
  If one is not in this extreme possibility, then 
  there would have to be 
 a broad and as yet unknown new class of
number fields which is also exceptional
from the point of view of the mass heuristic.

  \end{document}